\documentstyle[10pt, amscd]{amsart} 
 \renewcommand{\epsilon}{\varepsilon}
 \newcommand{\newsection}[1]
 {\subsection{#1}\setcounter{theorem}{0} \setcounter{equation}{0}
 \par\noindent}

 \newtheorem{theorem}{Theorem}
 
 \newtheorem{lemma}[theorem]{Lemma}
 \newtheorem{corr}[theorem]{Corollary}
 
 \newtheorem{proposition}[theorem]{Proposition}
 \newtheorem{deff}[theorem]{Definition}
 
 \newcommand{\bth}{\begin{theorem}}
 \newcommand{\ble}{\begin{lemma}}
 \newcommand{\bcor}{\begin{corr}}
 \newcommand{\bdeff}{\begin{deff}}
 \newcommand{\bprop}{\begin{proposition}}
 \newcommand{\eth}{\end{theorem}}
 \newcommand{\ele}{\end{lemma}}
 \newcommand{\ecor}{\end{corr}}
 \newcommand{\edeff}{\end{deff}}
 
 \newcommand{\eprop}{\end{proposition}}

 \renewcommand{\Pi}{\varPi}

 \renewcommand{\epsilon}{\varepsilon}

 \newcommand{\Rt}{{\Bbb R}^3}

 \newcommand{\Kob}{{\cal K}}
 \newcommand{\Dia}{\overline{\Bbb E}^{1+3}} \newcommand{\Diap}{\overline{\Bbb E}^{1+3}_+} \newcommand{\Cyl}{{\Bbb E}^{1+3}} \newcommand{\Cylp}{{\Bbb E}^{1+3}_+}
 \newcommand{\Penrose}{{\cal P}}
 \newcommand{\Rplus}{{\Bbb R}_+}
 
 \newcommand{\tidle}{\tilde}

\begin{document}

 \title[Global existence for nonlinear wave equations]{On Global existence for nonlinear
 wave equations outside of convex obstacles}
 \thanks{The authors were supported in part by the NSF}
 \author{Markus Keel}
 \author{Hart F. Smith}
 \author{Christopher D. Sogge}
 \address{Department of Mathematics, California Institute of Technology, Pasadena, CA 91125}
 \address{Department of Mathematics, University of Washington, Seattle, WA 98195}
 \address{Department of Mathematics, The Johns Hopkins University, Baltimore,
 MD 21218}
 
 \maketitle

 \newsection{Introduction}
 
 The purpose of this paper is to prove global existence for certain small-amplitude
 nonlinear Dirichlet-wave equations
 outside of smooth, strictly convex obstacles $\Kob \subset \Rt$.  As in earlier works on the boundaryless case studied by Christodoulou \cite{C}
 and Klainerman \cite{K} we shall be concerned with equations where the nonlinearities involve a null form.
 
 The null forms that we shall consider are the standard ones, which
 take the form
 \begin{equation}\label{Q0}
 Q_0(dv,dw)=\partial_tv\partial_tw-\sum_{j=1}^3\partial_jv\partial_jw
 \end{equation}
 or
 \begin{equation}\label{Qjk}
 Q_{jk}(dv,dw)=\partial_jv\partial_kw-\partial_kv\partial_jw, \quad 0\le j<k\le 3.
 \end{equation}
 Here, $\partial_jv=\partial v/\partial x_j$ if $j=1,2,3$, and $\partial_0v=\partial_tv=\partial v/\partial t$.
 Recall that nonlinear hyperbolic systems that satisfy Klainerman's \cite{K0} null condition must involve
 linear combinations of these null forms.  Furthermore, it is known (see \cite{C} and \cite{K}) that in
 Minkowski space there is global existence for small compactly supported data, while for other types of
 nonlinearities one can find arbitrarily small data with fixed compact support for which there is blowup. (See \cite{J}.)
 
 Based on this we shall study nonlinear hyperbolic equations in the exterior of a strictly convex domain
 in $\Rt$ that satisfy the null condition.  To be more specific, we shall fix a strictly convex obstacle
 $\Kob \subset \Rt$ with smooth boundary $\partial\Kob$.

 We then shall consider nonlinear systems of the form
 \begin{equation}\label{1.3}
 \begin{cases}
 \square u=Q(du,du), \quad (t,x)\in {\Bbb R}_+\times \Rt\backslash \Kob
 \\
 u(t,\cdot)|_{\partial\Kob}=0
 \\
 u(0,\cdot)=f, \, \, \partial_tu(0,\cdot)=g,
 \end{cases}
 \end{equation}
 if $\square = \partial^2/\partial t^2-\Delta$ is the usual D'Alembertian.
 Here $u=(u^1,u^2,\dots,u^N)$, and $Q=(Q^1,\dots,Q^N)$ with the $i$-th component of $Q$ being
 a linear combination of the null forms in \eqref{Q0} and \eqref{Qjk}.  More specifically,
 \begin{equation}\label{Qdef}
 Q^i=\sum_{j,k}a^i_{j,k}B^i_{j,k}(du^j,du^k),
 \end{equation}
 with $a^i_{j,k}$ being a constant and $B^i_{j,k}$ being any one of the null forms in
 \eqref{Q0}-\eqref{Qjk}.
 
 We shall be interested in obtaining $H^2$ solutions of our equation.  In view
 of the Dirichlet condition in \eqref{1.3} one must also then have that
 $\partial_tu(t,\cdot)|_{\partial\Kob}=0$. As a consequence, one must assume
 that the data satisfy the compatibility conditions
 \begin{equation}\label{compatibility}
 f(x)=0 \quad \text{and} \quad g(x)=0 \quad \text{if } \, x\in \partial\Kob.
 \end{equation}
 
 As in Christodoulou's \cite{C} results for the non-obstacle case, we shall not have to assume
 that the data has compact support.  Instead, this condition is replaced by a natural assumption
 that $f$ and $g$ belong to certain weighted Sobolev spaces.
 
 To be more precise, let us first recall the weighted Sobolev spaces that were
 used by Christodoulou \cite{C}, which are given by the norm
 $$\|f\|_{H^{m,j}(\Rt)}=\sum_{|\alpha|\le
 m}(\int_{\Rt}(1+|x|^2)^{|\alpha|+j}|D^\alpha f(x)|^2 \, dx )^{1/2}.$$ Here we
 are using the standard notation $D^\alpha=(\partial/\partial x_1)^{\alpha_1}
 (\partial/\partial x_2)^{\alpha_2}(\partial/\partial x_3)^{\alpha_3}$, and
 $|\alpha|=\alpha_1+\alpha_2+\alpha_3$.  (Later in this paper $D^\alpha$ will involve time
 derivatives, but this will be clear from the context.)  The associated weighted
 Dirichlet-Sobolev spaces for $m=1,2\dots$ then are
 \begin{equation}\label{euclidsobolev}
 H^{m,j}_D(\Rt\backslash \Kob)=\{f\in H^{m,j}(\Rt\backslash \Kob): \, f|_{\partial\Kob}=0\},
 \end{equation}
 where $H^{m,j}(\Rt\backslash \Kob)$ is the space of restrictions of  elements of $H^{m,j}(\Rt)$.  Hence,
 \begin{equation}\label{sobolev}
 \|f\|^2_{H^{m,j}_D}=\sum_{|\alpha|\le m}\int_{\Rt\backslash \Kob} (1+|x|^2)^{|\alpha|+j}|D^\alpha f(x)|^2 \, dx,
 \end{equation}
gives the natural norm on $H^{m,j}_D(\Rt\backslash \Kob)$. Also,
 $H^m({\Bbb R}^3\backslash \Kob)$, $m=1,2,3,\dots$ will denote the Sobolev space
 of restrictions of elements of $H^m({\Bbb R}^3)$, while $H^m_D({\Bbb
 R}^3\backslash \Kob)$, $m=1,2,3,\dots$ will denote the subset $\{f\in H^m({\Bbb
 R}^3\backslash \Kob): f|_{\partial \Kob}=0\}$.
 
 We can now state our first main result.
 
 \begin{theorem}\label{theorem1.1}  Suppose that $f$ and $g$ satisfy \eqref{compatibility} and that
 \begin{equation}\label{1.8}
 \|f\|_{H^{2,1}_D}+\|g\|_{H^{1,2}_D}<\varepsilon_0.
 \end{equation}
 Then if $\varepsilon_0>0$ is sufficiently small there is a unique global solution $u$ of \eqref{1.3} verifying
 \begin{equation}\label{1.9}
 \sup_{0<t<T}\sum_{|\alpha|\le2}\|D^\alpha u(t,\cdot)\|_{L^2(\Rt\backslash \Kob)}
 +\sum_{|\alpha|\le1}\|D^\alpha Q\|_{L^2([0,T]\times \Rt\backslash \Kob)}<\infty, \quad T>0.
 \end{equation}
 Furthermore, the solution has the decay property
 \begin{equation}\label{1.10}
 |u(t,x)|\le C(1+t)^{-1}.
 \end{equation}
 \end{theorem}
 
 %We shall also see that we can get global in time
 %smooth solutions if the data is smooth and satisfies the
 %appropriate compatibility conditions.

 Analogous results under the assumption of spherical symmetry for $\Kob$ and $u$ were obtained
 by Godin \cite{G} if $Q=Q_0$.  His proof involved an adaptation of Christodoulou's \cite{C}
 method to this setting.  If one drops the assumption of spherical symmetry, it does not seem
 that the arguments in \cite{G} will apply in a straightforward way.
 
 Previous work in higher dimensions applied Lorentz vector field techniques
 to the exterior problem.  For general nonlinearities
 quadratic in $\nabla u, \partial_t u$, global smooth solutions were
 shown by Shibata and Tsutsumi \cite{ST} to exist in dimension $ n \geq 6$.  In 
Hayashi
\cite{Hanew}, global existence of smooth solutions outside of spheres in $n\ge 4$
is shown for a restricted class of quadratic nonlinearities, extending work
in \cite{Ha}.
%\cite{Ha}, global
% existence of smooth solutions in $n \geq 4$ is shown for a restricted
% class of quadratic nonlinearities under spherical symmetry assumptions
% on both data and obstacle.
 
 Our methods also give smooth solutions if the data satisfies the necessary
 compatibility conditions.  To state these conditions, we define a collection of
 functions $\psi_j$ on ${\Bbb R}^3\backslash \Kob$ as follows.  Set $$\psi_0=f,
 \quad \psi_1=g.$$ We now define $\psi_j$ recursively so that, if the function
 $u_c$ has the following Taylor expansion in $t$, $$u_c(t,x)\approx
 \sum_{j=0}^\infty\psi_j(x)t^j/j!,$$ then $\square u_c-Q(du_c,du_c)$ vanishes to
 infinite order at $t=0$. This is seen to determine $\psi_j$ uniquely, and
 $\psi_j$ is a nonlinear function of the data $(f,g)$ that involves derivatives
 of order up to $j$ of $f$, and of order up to $j-1$ of $g$.  If $(f,g)\in
 H^k({\Bbb R}^3\backslash \Kob)\times H^{k-1}({\Bbb R}^3\backslash \Kob)$ one
 would only use the first $k$ terms and require that the resulting function
 vanishes to order $k$ at $t=0$ if one wishes to obtain $H^k$ solutions.
 
 \begin{deff}\label{compatibilitydef}  We say that the data $(f,g)$
 satisfies the compatibility conditions to order $k$ if, for $0\le j\le k$, the
 functions $\psi_j$ vanish on $\partial \Kob$.  We say that the data satisfies
 the compatibility conditions to infinite order if all the $\psi_j$ vanish on
 $\partial\Kob$.
 \end{deff}
 
 Note that the assumption in Theorem \ref{theorem1.1} was that the data
 satisfies the first order compatibility condition. If $\overline{f}$ and
 $\overline{g}$ are appropriate extensions of $f$ and $g$, respectively, to all
 of ${\Bbb R}^3$ and if $\overline{u}$ is a local solution to the equation
 $\square \overline{u}=Q(d\overline{u},d\overline{u})$ one can just take
 $\psi_j(x)=\partial^j_t \overline{u}(0,x)$, $x\in {\Bbb R}^3\backslash \Kob$.
 Then the $k$th order compatibility condition would be equivalent to the
 condition that $\partial_t^j\overline{u}(0,x)=0$, when $0\le j\le k$ and $x\in
 \partial\Kob$. Also note that all of the compatibility conditions are
 automatically satisfied if the data vanishes near $\partial\Kob$.
 
 \begin{theorem}\label{theorem1.3}  Assume that $(f,g)\in
 C^\infty({\Bbb R}^3\backslash \Kob)$ satisfies the compatibility conditions to
 infinite order and that
 \begin{equation}\label{highersmall}
 \|f\|_{H^{2,1}_D}+\|g\|_{H^{1,2}_D}<\varepsilon_0.\end{equation} Then if
 $\varepsilon_0>0$ is sufficiently small there is a unique solution $u\in
 C^\infty({\Bbb R}_+\times {\Bbb R}^3\backslash \Kob)$ of \eqref{1.3}.  This
 solution also satisfies \eqref{1.10}.\end{theorem}

 Our approach will be to combine the conformal method introduced in \cite{C1}
 (see also \cite{C})
with techniques developed in \cite{GS}, \cite{KM}, \cite{SS1} and \cite{S1}.
 These results all involve generalizations of an inequality of Klainerman and
 Machedon \cite{KM}.  In the scalar case, it involves estimates for solutions
 $u$ of linear wave equations in Minkowski space ${\Bbb R}^{1+3}_+$
 \begin{equation}\begin{cases}\label{minkinhom}
 \square u=F, \quad (t,x)\in {\Bbb R}_+^{1+3}
 \\
 u(t,\cdot)=u_0, \quad \partial_tu(t,\cdot)=u_1.
 \end{cases}
 \end{equation}
 To be more specific, if $v$ solves the same equation with data $(v_0,v_1,G)\,,$ then the main
 estimate in \cite{KM} says that
 if $Q$ is any of the null forms in \eqref{Q0}-\eqref{Qjk}
 then
 \begin{multline}\label{KM}
 \sum_{|\alpha|\le1}\|D^\alpha Q(du,dv)\|_{L^2({\Bbb R}^{1+3}_+)}
 \le
 C\bigl(\|u_0\|_{H^2}+\|u_1\|_{H^1}+\int_0^\infty \|F(t,\cdot)\|_{H^1}\,
 dt\bigr)
 \\
 \times \bigl(\|v_0\|_{H^2}+\|v_1\|_{H^1}+\int_0^\infty
 \|G(t,\cdot)\|_{H^1}\, dt\bigr).
 \end{multline}
 
 This estimate was extended locally to the boundaryless manifold case by the last author in \cite{S1}.
 Since we shall use a variant of Christodoulou's conformal method,
 this will be one of the ingredients in the proof of our existence results.  In particular, the results from \cite{S1}
 in the special case where the underlying spatial manifold is $S^3$ will allow us to prove the
 necessary estimates when we are away from the boundary in the image of
 ${\Bbb R}_+ \times {\Bbb R}^3\backslash \Kob$
 via Penrose's \cite{P} conformal compactification.

 Proving estimates near the boundary in the image of
 ${\Bbb R}_+ \times {\Bbb R}^3\backslash \Kob$, though, is much more complicated.  This is
 because of the fact that the fixed-time cross sections of this set are
 hypersurfaces that vary with time in $[0,\pi)\times S^3$, and, in fact, degenerate to a point
 as the time variable $T$ tends to $\pi$.
 
 To get around this it turns out to be convenient to pull back our estimates for the image of
 ${\Bbb R}_+\times {\Bbb R}^3_+\backslash \Kob$ to Minkowski space.  Here we shall be able to
 prove the resulting estimates by applying energy decay estimates of Lax, Morawetz and
 Phillips \cite{LMP}-\cite{LP} along with local estimates for the inhomogeneous Dirichlet-wave
 equation.  The latter were recently obtained by the last two authors in \cite{SS2}.
 
 Let us be more specific.
 
 The classical energy estimates from \cite{LMP}-\cite{LP} that we shall require are contained in the following
 
 \begin{theorem}\label{energydecaytheorem}  (Lax-Morawetz-Phillips \cite{LMP}-\cite{LP})  Suppose
 that $u$ solves the homogeneous Dirichlet-wave equation $\square u=0$ in
 ${\Bbb R}_+\times \Rt\backslash \Kob$, $u|_{{\Bbb R}_+\times \partial
 \Kob}=0$.  Suppose further that the Cauchy data $u_0=u|_{t=0}$ and $u_1=\partial_t u|_{t=0}$ vanish
 when $|x|>A$, with $A$ fixed.  Then there are constants $c>0$ and $C<+\infty$, depending on $\Kob$
 and $A$, so that if $Du=(\partial_tu,\nabla_xu)$
 \begin{equation}\label{1.13}
 \sum_{|\alpha|\le 1}\int_{\{x\in \Rt\backslash \Kob: |x|<A\}}|D^\alpha u(t,x)|^2\, dx
 \le Ce^{-ct}\sum_{|\alpha|\le1}
\int_{\{x\in \Rt\backslash \Kob: |x|<A\}}|D^\alpha u(0,x)|^2\, dx.
 \end{equation}
 \end{theorem}
 
 We shall actually require a consequence of \eqref{1.13} that involves estimates for higher
 derivatives.  These follow from \eqref{1.13} and elliptic regularity arguments (cf.
 \cite{ST}).  If we use Duhamel's formula we can also get estimates for the inhomogeneous
 wave equation
 \begin{equation}\label{1.14}
 \begin{cases}
 \square u(t,x)=F(t,x), \quad (t,x)\in {\Bbb R}_+\times {\Bbb R}^3\backslash \Kob
 \\
 u(t,x)=0, \quad x\in \partial \Kob
 \\
 u(0,x)=u_0(x), \quad \partial_t u(0,x)=u_1(x).
 \end{cases}
 \end{equation}
 Specifically, we first notice that \eqref{1.13} yields
 \begin{multline*}
 \sum_{|\alpha|\le2}\|D^\alpha u(t,\cdot)\|_{L^2(\{x\in {\Bbb R}^3\backslash
 \Kob : |x|<A\})}
 \\
 \le Ce^{-ct}\bigl[\|u_0\|_{H^2_D}+\|u_1\|_{H^1_D}+\sum_{|\alpha|\le1}\|e^{cs}D^\alpha
 F\|_{L^1_tL^2_x([0,t]\times {\Bbb R}^3\backslash \Kob)}\bigr]
 ,
 \end{multline*}
 for constants $c>0$ and $C<+\infty$ as above
 provided that $u_j(x)=0$, $j=0,1$, and $F(t,x)=0$ when $|x|>A$.  Note that if we use the
 Schwarz inequality then we can dominate the last term by
 $$C\sum_{|\alpha|\le1}\|e^{-c(t-s)/2}D^\alpha F\|_{L^2([0,t]\times {\Bbb R}^3\backslash
 \Kob)}.$$
 If we combine the last two inequalities we obtain the following useful result.
 
 \begin{proposition}\label{prop1.3}  Fix $A$ and suppose that $u$ solves \eqref{1.14} with Cauchy
 data $u_j(x)$, $j=0,1$, and forcing term $F(t,x)$ both vanishing when $|x|>A$.  Then there
 are constants $c>0$ and $C<+\infty$ depending only on $A$ and $\Kob$ so that
 \begin{multline}\label{energydecay}
 \sum_{|\alpha|\le2}\|D^\alpha u(t,\cdot)\|_{L^2(\{x\in {\Bbb R}^3\backslash \Kob: |x|<A\})}
 \\
 \le Ce^{-ct}\bigl(
 \|u_0\|_{H^2_D}+\|u_1\|_{H^1_D}\bigr)+C\sum_{|\alpha|\le1}
 \|e^{-c(t-s)}D^\alpha F\|_{L^2([0,t]\times {\Bbb R}^3\backslash \Kob)}.
 \end{multline}
 \end{proposition}
 
 The local null form estimates that we shall need are:
 
 \begin{theorem}\label{theorem1.6}
 (Smith-Sogge \cite{SS2})  Suppose that $u$ and $v$ satisfy inhomogeneous Dirichlet-wave equations
 $\square u=F$ and $\square v=G$ in ${\Bbb R}_+\times \Rt\backslash \Kob$ with Cauchy data $(u_0,u_1)$
 and $(v_0,v_1)$, respectively.  Then if $Q$ is any of the null forms in \eqref{Q0}-\eqref{Qjk}
 \begin{multline}\label{1.17}
 \sum_{|\alpha|\le1}\|D^\alpha Q(du,dv)\|_{L^2([0,1]\times \Rt\backslash \Kob)}
 \\
 \le C
 \Bigl(
 \|u_0\|_{H^2_D}+\|u_1\|_{H^1_D}+
 \sum_{|\alpha|\le1}\|D^\alpha F(t,\cdot)
 \|_{L^1L^2([0,1]\times \Rt\backslash \Kob)}\Bigr)
 \\  \times
 \Bigl(
 \|v_0\|_{H^2_D}+\|v_1\|_{H^1_D}+
 \sum_{|\alpha|\le1}\|D^\alpha G(t,\cdot)
 \|_{L^1L^2([0,1]\times \Rt\backslash \Kob)}\Bigr).
 \end{multline}
 \end{theorem}
 
 This equation of course gives small $H^2$ data local existence for \eqref{1.3} (see
 \cite{SS2}.)  Unfortunately, though, it does not yield global existence. To get around this, we shall prove a weighted variant in the
 image of ${\Bbb R}_+\times \Rt\backslash \Kob$ in the Einstein diamond. (This estimate would in turn pull back to a global weighted
 variant of \eqref{1.17}, but we do not explore that in this paper.)  
 This paper is organized as follows.  In the next section we shall review
 Penrose's \cite{P} conformal compactification of Minkowski space to the so
 called Einstein diamond in $(-\pi,\pi)\times S^3$. We shall collect the necessary
 facts regarding the way that our nonlinear equation \eqref{1.3} transforms and
 state our main estimate in $[0,\pi)\times S^3$ that leads to global existence.
 It will be a weighted analog of \eqref{1.17} on $[0,\pi)\times S^3$
 minus the image of ${\Bbb R}_+\times \Kob$. The weights will compensate for the
 degeneracies of the boundary as $T\to \pi$.  In subsequent sections we shall
 prove our main estimate using the strategy mentioned above of proving things
 directly outside of an appropriate neighborhood of the boundary, while proving
 things near the boundary by pulling everything back to Minkowski space.
 Finally, after we prove the weighted estimates we shall see how our estimates
 give the $H^2$ and $C^\infty$ global existence theorems.
 
 \newsection{Conformal Compactification and the Main Estimate}
 
 Consider polar coordinates on the sphere $S^3$ minus the south pole given by
 \begin{equation}\label{polar}
 [0,\pi)\times S^2\ni (R,\omega)\to (\cos R, \omega \sin R) = X = (X_0, \overrightarrow{X})\in {\Bbb R}^4.
 \end{equation}
 Then the ``standard'' Lorentz metric on ${\Bbb R}\times S^3$ is given by
 \begin{equation}\label{standardmetric}
 g=dX^2=dT^2-dR^2-\sin^2R\, d\omega^2.
 \end{equation}
 
 Since we are interested in solutions of \eqref{1.3} it is natural to consider another metric coming from Penrose's
 conformal compactification of ${\Bbb R}^{1+3}={\Bbb R}\times \Rt$.  This ``physical metric'' is just the pushforward
 of the standard Lorentz metric in Minkowski space, and it turns out to be conformally equivalent to \eqref{standardmetric}.
 The facts we shall state about this transformation are well known and can be found, for example, in H\"ormander \cite{H}.
 
 Let us be more specific.  First of all, we define the Einstein cylinder to be the set $$ \Cyl=(-\pi,\pi)\times S^3\,. $$ We also need to define the Einstein diamond
 $$
 \Dia \subset \Cyl
 $$
 which is just the proper subset of $(-\pi,\pi)\times S^3$ given by
 \begin{equation}\label{einstein}
 \Dia=\{(T,\cos R,\omega \sin R): \, -\pi < T<\pi, \, \, 0\le R< \pi, \, \, R+|T|<\pi\}.
 \end{equation}
 Penrose's transformation in \cite{P} then is the conformal map ${\cal P}:{\Bbb R}^{1+3}\to \Dia$,
 which in polar coordinates is defined by
 \begin{multline}\label{penroseradial}
 (t,r,\omega)\to (T,R,\omega),
 \\ T=\tan^{-1}(t+r)+\tan^{-1}(t-r), \, \, R=\tan^{-1}(t+r)-\tan^{-1}(t-r),
 \end{multline}
 where as usual $x=r\omega$, $r=|x|$ in ${\Bbb R}^3$.
 Note that the inverse of
 $\Penrose$ is given by
 \begin{equation}\label{inversepenrosetransform}
 {\cal P}^{-1}: \Dia\ni (T,X)\to \frac{1}{\cos T+X_0}(\sin T, \overrightarrow{X})\in {\Bbb R}^{1+3}.
 \end{equation}
 Notice also that when $T=0$ the map \eqref{inversepenrosetransform} is
 stereographic projection of $S^3$ from the south pole $(-1,0,0,0)$.
 
 Under this map the pushforward of the Minkowski metric $dt^2-dx^2$ is the Lorentz metric $\tilde g$ in
 $\Dia$ defined by
 \begin{equation}\label{metrics}
 g=\Omega^2\tilde g,
 \end{equation}
 where the conformal factor is given by
 \begin{equation}\label{conformal}
 \Omega=\cos T+\cos R=\frac2{(1+(t+r)^2)^{1/2}(1+(t-r)^2)^{1/2}},
 \end{equation}
 with $(t,r)$ and $(T,R)$ being identified as in \eqref{penroseradial}.
 
 Continuing, let
 $$\square_g=\partial^2/\partial T^2-\Delta_{S^3}$$
 be the D'Alembertian coming from the standard Laplace-Beltrami operator, $\Delta_{S^3}$,
 given by the round metric $dT^2-dX^2$.  On the other hand, $\tilde \square$ will be the D'Alembertian on
 ${\Bbb R}^{1+3}$ or $\Dia$, depending on the context, that arises from $dt^2-dx^2$.  In Minkowski
 space it is of course the standard D'Alembertian $\square =\partial^2/\partial t^2-\Delta$.
 A key fact for us will be the way the two D'Alembertians are related in $\Dia$:
 \begin{equation}\label{intertwine}
 \square_g+1=\Omega^{-3}\tilde \square \Omega,
 \end{equation}
 with the additive constant $1$ arising from the non-zero scalar curvature of $g$.
 Equivalently,
 \begin{equation}\label{key}
 (\square_g+1)u=F \, \, \iff \, \, \tilde \square \tilde u=\tilde F \, \,
 \text{with } \, \, \tilde u=\Omega u \, \, \text{ and } \, \, \tilde F=\Omega^3F.
 \end{equation}
 
 Another related standard fact concerns the way that the Sobolev spaces
 \eqref{euclidsobolev} transform.  For this we need to let ${\cal P}_0(\Kob)$
 denote the image of $\Kob$ under the restriction of the Penrose transformation
 to $\{(0,x):\, x\in \Rt\}$.  Since ${\cal P}_0$ is the stereographic projection
 map, ${\cal P}_0(\Kob)\subset S^3$ is a convex set with smooth boundary.  For
 $m=1,2,\dots$ we then let $H^m_D(S^3\backslash {\cal P}_0(\Kob)) =\{f\in
 H^m(S^3\backslash {\cal P}_0(\Kob)): \, f|_{\partial {\cal P}_0(\Kob)}=0\}$ be
 the Sobolev space defined as before.  If then ${\cal P}_0^*h$ denotes the
 pullback of a function $h$ on $S^3$ via the above restriction of the Penrose
 transformation it follows that the mapping
 \begin{equation*}\tilde f=\Omega {\cal P}_0^*f  \to f
 \end{equation*}
 is a continuous map from $H^{m,m-1}({\Bbb R}^3\backslash \Kob)$ to
 $H^m(S^3\backslash {\cal P}_0(\Kob))$.  That is, for every $m=1,2,3,\dots$
 there is a constant $C_m<\infty$ so that
 \begin{equation}\label{2.10}
 \|f\|_{H^m(S^3\backslash {\cal P}_0(\Kob))}\le C_m\|\Omega {\cal
 P}_0^*f\|_{H^{m,m-1}({\Bbb R}^3\backslash \Kob)}=C_m\|\tilde
 f\|_{H^{m,m-1}({\Bbb R}^3\backslash \Kob)}.
 \end{equation}
 Notice that if $u$ and $\tilde u$ are related as above then this inequality
 yields
 \begin{equation*}
 \|\partial_T u(0,\cdot)\|_{H^{m-1}(S^3\backslash {\cal P}_0(\Kob))} \le
 C_m\|\partial_t \tilde u(0,\cdot)\|_{H^{m-1,m}({\Bbb R}^3\backslash \Kob)},
 \end{equation*}
 since the pushforward of $\partial/\partial t$ is $\Omega \,\partial/\partial T$
 when $t=0$ and since $\Omega=2/(1+|x|^2)$ when $t=0$.
 %induces
 %an isomorphism of $H^m_D(S^3\backslash {\cal P}_0(\Kob))$ and $H^{m,m-1}_D(\Rt\backslash \Kob)$.  That
 %is, for every $m=1,2,\dots$ there is a constant $C=C_m\in (0,+\infty)$ so that
 %\begin{equation}\label{equivnorms}
 %C^{-1}\|f\|_{H^m_D(S^3\backslash {\cal P}_0(\Kob))}
 %\le \|\Omega {\cal P}_0^*f\|_{H^{m,m-1}_D(\Rt\backslash \Kob)}
 %\le C\|f\|_{H^m_D(S^3\backslash {\cal P}_0(\Kob))}.
 %\end{equation}
 
 Another fact that we shall use concerns the image of ${\Bbb R}_+\times {\cal K}$ in
 $\Diap=\Dia\cap [0,\pi)\times S^3$.  Let us call this set
 \begin{equation}\label{transformedboundary}
 \Kob_*=\Penrose(\Rplus \times \Kob).
 \end{equation}
 Then one can check using \eqref{inversepenrosetransform} that there is a uniform constant
 $0<C<\infty$ so that for $0\le T<\pi$
 \begin{equation}\label{boundary1}
 C^{-1}(\pi-T)^2\le \text{dist}(X,{\bold 1})\le C(\pi-T)^2 \quad \text{if } \, \, (T,X)\in \partial \Kob_*.
 \end{equation}
 Here
 $${\bold 1}=(1,0,0,0)$$
 is the north pole and $\text{dist}(\cdot\,, \cdot)$ is the standard distance on $S^3$ (induced by
 the metric $g$).  If we let
 \begin{equation}\label{P0}
 P_0=(\pi,{\bold 1})
 \end{equation}
 be the convergence of the positive time-like infinities then we can state \eqref{boundary1} in
 an equivalent way by saying that
 \begin{equation}\label{boundary2}
 C^{-1}\text{dist}^2(P,P_0)\le \text{dist}(X,{\bold 1})\le C\,\text{dist}^2(P,P_0),
 \, \text{if }  P=(T,X)\in \partial \Kob_*.
 \end{equation}
 
 It should be clear from the context that we are using $\text{dist}(\cdot\,, \cdot)$ in two different ways.
 The distance between two points on $S^3$ is the distance given by the restriction of $g$ to $S^3$ and the
 distance between two points on ${\Bbb R}\times S^3$ is given by $g$.  We shall use this notation in what follows.
 
 The fact that the boundary of $\Kob*$ varies with time and, moreover, degenerates to a point as
 $T\to \pi$ is the reason that Christodoulou's \cite{C} approach of using the above conformal
 compactification and the energy integral method does not seem to apply in an easy way for
 \eqref{1.3}.  In particular, as we shall see, it does not seem easy to control fixed-time high order
 Sobolev norms of a solution of the Dirichlet-wave equation in $\Diap\backslash \Kob_*$
 as $T\to \pi$, unless one is willing to include weights in the norms involving powers of $\text{dist}^2(P,P_0)$.
 
 We are almost ready to state our main inequality.  To motivate the weights we shall use, let us first recall
 how the standard vector fields in Minkowski space pushforward to ones in $\Dia$ via ${\cal P}$.
 We shall follow the exposition in H\"ormander \cite{H}, p. 277-282.
 
 To do this it is convenient to use stereographic projection coordinates on $S^3$.  As we pointed out before the
 south pole stereographic projection coordinates come from the restriction of $\Penrose^{-1}$ to $T=0$:
 \begin{equation}\label{southcoord}
 Y=\Penrose_0^{-1}(X)=\frac{\sin R}{1+\cos R}\omega=\tan (\tfrac{R}{2}) \omega.
 \end{equation}
 The coordinates coming from the stereographic north pole projection arise from these and the Kelvin transform
 \begin{equation}\label{northcord}
 Z_j=Y_j/|Y|^2.
 \end{equation}
 
 To compute the pushforwards of vector fields on $\Dia$ it is convenient to use the vector fields
 \begin{equation}\label{spherevectorfields}
 \frac{\partial}{\partial T}, \, \, X_j\frac{\partial}{\partial X_k}-X_k\frac{\partial}{\partial X_j}, \, \, \,
 0\le j<k\le 3.
 \end{equation}
 For future reference, let us arrange these as
 $\Gamma=\{\Gamma_0,\dots,\Gamma_6\}$
 and write $\Gamma^\alpha = \Gamma_0^{\alpha_0}\cdots\Gamma_6^{\alpha_6}$. Note that each $\Gamma$ is actually a smooth vector field on $\Cyl$.

 With this notation we can state the following result (see \cite{H})
 
 \bprop\label{derivativesprop}
 The pushforwards of $\partial/\partial t$ and $\partial/\partial x_j$, $1\le j\le 3$ in ${\Bbb R}^{1+3}$
 by ${\cal P}$ are
 \begin{align}\label{dt}
 \frac\partial{\partial t}= &\bigl(1+\frac{1-|Y|^2}{1+|Y|^2}\cos T\bigr)\tfrac\partial{\partial T}
 -\sin T\langle Y,\tfrac\partial{\partial Y}\rangle
 \\
 \label{rrrr}
 =&
 \bigl(1+\frac{|Z|^2-1}{|Z|^2+1}\cos T\bigr)\tfrac\partial{\partial T}
 +\sin T\langle Z,\tfrac\partial{\partial Z}\rangle,
 \end{align}
 and
 \begin{align}\label{dx}
 \frac{\partial}{\partial x_j}
 &=\frac{-2Y_j}{1+|Y|^2}\sin T\tfrac{\partial}{\partial T}+\tfrac12((1+|Y|^2)\cos T+1-|Y|^2)\frac{\partial}{\partial Y_j}
 +(1-\cos T)Y_j\langle Y,\frac{\partial}{\partial Y}\rangle
 \\
 \label{zzz}
 &=\frac{-2Z_j}{1+|Z|^2}\sin T\tfrac{\partial}{\partial T}+\tfrac12((1+|Z|^2)\cos T+|Z|^2-1)\tfrac\partial{\partial Z_j}
 +(1+\cos T)Z_j\langle Z,\tfrac\partial{\partial Z}\rangle.
 \end{align}
 The pushforward of the vector fields in \eqref{spherevectorfields} via ${\cal P}^{-1}$ are given by
 \begin{align}\label{Xjk}
 X_j\frac\partial{\partial X_k}-X_k\frac{\partial}{\partial X_j}&=x_j\frac{\partial}{\partial x_k}-
 x_k\frac{\partial}{\partial x_j}, \, \, 1\le j<k\le 3,
 \\ \label{X0k}
 X_0\frac{\partial}{\partial X_k}-X_k\frac{\partial}{\partial X_0}
 &=\tfrac12(1+t^2-|x|^2)\frac{\partial}{\partial x_k}+
 x_k(t\frac{\partial}{\partial t}+\langle x,\frac{\partial}{\partial x}\rangle), \, \, 1\le k\le 3,
 \end{align}
 and
 \begin{equation}\label{dT}
 \frac{\partial}{\partial T}=\tfrac12(1+t^2-|x|^2)\frac\partial{\partial t}+t\langle x,\frac\partial{\partial x}\rangle.
 \end{equation}
 %Check equation dT
 
 \eprop
 
 Note that the coefficients of $\partial/\partial T$ and $\partial/\partial Y$ in \eqref{dt} and
 \eqref{dx} are $O((\pi-T)^2)$ when, say, $0\le T<\pi$ and $R\le (\pi-T)/4$.  Similarly, if $|x|\le t/4$ the
 coefficients of $\partial/\partial t$ and $\partial/\partial x$ in \eqref{X0k} and \eqref{dT}
 are $O(t^2+|x|^2)$.  Hence we have the following useful result:
 
 \begin{proposition}\label{middlebehavior}  In the region where $|x|\le t/4$ we can write
 $$\frac\partial{\partial t}=\sum a_{0k}(T,X)\Gamma_k,
 \quad \text{and }\, \,  \, \,
 \frac{\partial}{\partial x_j}=\sum a_{jk}(T,X)\Gamma_k,$$
 where if $P_0$ is as in \eqref{P0}
 $$|\Gamma^\alpha a_{jk}|\le C\,\text{dist}((T,X),P_0)^{2-|\alpha|}, \quad |\alpha|\le2.
 $$
 Also, if $0\le T<\pi$ and $R\le (\pi-T)/4$
 $$\frac{\partial}{\partial T}=b_{00}(t,x)\frac{\partial}{\partial t}+\sum b_{0k}(t,x)\frac{\partial}{\partial x_k}$$
 and
 $$X_0\frac{\partial}{\partial X_j}-X_j\frac{\partial}{\partial X_0}=
 b_{0j}(t,x)\frac{\partial}{\partial t}+\sum b_{jk}(t,x)\frac{\partial}{\partial x_j},$$
 where if $D=(\partial/\partial t,\partial/\partial x_1,\dots,\partial/\partial x_3)$
 $$|D^\alpha b_{jk}|\le C(1+|(t,x)|)^{2-|\alpha|}, \quad |\alpha|\le2.
 $$
 \end{proposition}

 We need one more thing before we can state our main inequality.  If $u$ is a function on  $\Diap$,
 as before, let $\tilde u$ denote the pullback of $\Omega u$ to ${\Bbb R}^{1+3}$.  If we fix a null form $Q$
 as in \eqref{Q0} or \eqref{Qjk} we shall let
 \begin{multline}\label{calQ}
 {\cal Q}(u(T,X), du(T,X);v(T,X),dv(T,X))
 \\
 =\Omega^{-3}Q(d\tilde u(t,x),d\tilde
 v(t,x)), \quad \Penrose(t,x)=(T,X).
 \end{multline}
 The reason for this is that if we change the notation a bit from the preceding section
 and write our main equation \eqref{1.3} as
 \begin{equation}\label{tildeeqn}
 \begin{cases}
 \tilde \square \tilde u=Q(d\tidle u,d\tilde u), \quad (t,x)\in {\Bbb R}_+\times {\Bbb R}^3\backslash {\cal K}
 \\
 \tilde u(t,\cdot)|_{\partial \Kob}=0
 \\
 \tilde u(0,\cdot)=\tilde f, \quad \partial_t \tilde u(0,\cdot)=\tilde g,
 \end{cases}
 \end{equation}
 then by \eqref{key} it transforms to the following equation in $\Diap\backslash \Kob_*$
 \begin{equation}\label{einsteineqn}
 \begin{cases}
 (\square_g +1)u={\cal Q}(u, du; u, du)
 \\
 u(T,X)=0, \quad (T,X)\in \partial \Kob_*
 \\
 u(0,\cdot)=f, \quad \partial_Tu(0,\cdot)=g,
 \end{cases}
 \end{equation}
 if the data satisfies $\tilde f=\Omega f$ and $\tilde g=\Omega^2 g$ and if $\Kob_*$ is as in
 \eqref{transformedboundary}.

 We are now in a position to state our main estimate.  It can be thought of as
 an appropriate version of Theorem \ref{theorem1.6} for $\Cylp\backslash
 \Kob_*$, where $\Kob_*$ is as above. Also, we shall see in Section 7
 that it immediately yields the desired existence results in Theorem
 \ref{theorem1.1}.  It involves solutions of the Dirichlet-wave equation in
 $\Cylp \backslash \Kob_*$:
 \begin{equation}\label{linear}
 \begin{cases}
 (\square_g+1)u=F,
 \\
 u|_{\partial \Kob_*}=0
 \\
 u|_{T=0}=u_0, \quad \partial_Tu|_{T=0}=u_1.
 \end{cases}
 \end{equation}
 
 \begin{theorem}\label{theorem2.3}
 Suppose that $u$ solves \eqref{linear} and that $v$ solves the same equation with $u_0$, $u_1$ and $F$
 replaced by $v_0$, $v_1$ and $G$, respectively.  If $P_0$ is as in \eqref{P0}, then
 \begin{multline}\label{nullest}
 \sum_{|\alpha|\le1}\|(\text{dist}^2(P,P_0)\Gamma)^\alpha {\cal
 Q}(u,du;v,dv)\|_{L^2(\Cylp\backslash \Kob_*)}
 \\
 \le
 C\bigl(\|u_0\|_{H^2_D}+\|u_1\|_{H^1_D}
 +\sum_{|\alpha|\le1}\|(\text{dist}^2(P,P_0)\Gamma)^\alpha F\|_{L^2(\Cylp\backslash \Kob_*)}\bigr)
 \\
 \times
 \bigl(\|v_0\|_{H^2_D}+\|v_1\|_{H^1_D}
 +\sum_{|\alpha|\le1}\|(\text{dist}^2(P,P_0)\Gamma)^\alpha G\|_{L^2(\Cylp\backslash \Kob_*)}\bigr).
 \end{multline}
 Furthermore, the following estimate holds \begin{multline}\label{weightedenergy}
 \sum_{|\alpha|\le 1}
 \|\text{dist}^{2|\alpha|}(P,P_0)\Gamma^\alpha u(T,\cdot)\|_{L^8(\Cylp\backslash \Kob_*)}\\
 \le C\bigl(\|u_0\|_{H^2_D}+\|u_1\|_{H^1_D}+\sum_{|\alpha|\le1}\|(\text{dist}^2(P,P_0)\Gamma)^\alpha
 F\|_{L^2(\Cylp\backslash \Kob_*)}\bigr).
 \end{multline}
 \end{theorem}

 In many ways \eqref{nullest} is a natural extension of \eqref{1.17} to the
 current setting.  This is because, near $\partial \Kob_*$, Minkowski
 derivatives $\partial/\partial x_j$, $0\le j\le 3$ are pushed forward via
 $\Penrose$ to ones in $\Dia \backslash \Kob_*$ that behave like linear
 combinations of the $\Gamma$ with $O(\text{dist}^2(P,P_0))$ coefficients.  In
 view of \eqref{boundary2}, the weights in \eqref{nullest} also naturally
 compensate for the degeneracy of $\partial \Kob_*$ as $T\to \pi$.
 
 The proof of \eqref{nullest} will require more precise information about the behavior of ${\cal Q}$ as $P\to P_0$.
 The result we shall need says in part that ${\cal Q}$ involves the standard null forms on ${\Bbb R}\times S^3$.
 These, we recall (see \cite{C}, \cite{S1}, \cite{GS}) are
 $$\overline{Q}_0(du,dv)=\sum g^{\mu\nu}(T,X)u_\mu v_\nu,$$
 and
 $$\overline{Q}_{jk}(du,dv)=\Gamma_ju\Gamma_kv-\Gamma_ku\Gamma_jv,$$
 where $u_\mu$ denotes the differential of $u$, $g^{\mu\nu}$ is the cometric associated with $g$,
 and $\Gamma_j$, $0\le j\le 6$, are as above.
 
 With this notation, we have the following
 
 \begin{proposition}\label{nulltransformationprop}
 Fix a ${\Bbb R}^{1+3}$ null form $Q$ in \eqref{Q0}-\eqref{Qjk} and let ${\cal Q}(u,du;v,dv)$ be defined
 as in \eqref{calQ}.
 Then ${\cal Q}$ extends to a bilinear function of $(u,dv; v,dv)$ on ${\Bbb R}\times S^3$ with smooth
 coefficients.  Moreover, we can write
 \begin{equation}\label{Qformula}
 {\cal Q}= a_0\overline{Q}_0(du,dv)+\sum_{jk}a_{jk}\overline{Q}_{jk}(du,dv)
 +u \sum_{|\alpha|=1} b_{1,\alpha}\Gamma^\alpha v
 +v \sum_{|\alpha|=1} b_{2,\alpha}\Gamma^\alpha u
 +c(T,X) uv,
 \end{equation}
 where the coefficients are smooth, and, moreover,
 $$\Gamma^\alpha a_0, \, \, \Gamma^\alpha a_{jk} = O(\text{dist}^{2-|\alpha|}(P,P_0)),
 \quad 0\le |\alpha|\le 2,$$
 and
 $$b_{j,\alpha}=O(\text{dist}(P,P_0)).$$
 \end{proposition}
 
 The proof of Proposition \ref{nulltransformationprop} has two steps.  The hard step was carried out by Christodoulou \cite{C} who showed that
 one can write ${\cal Q}$ as in \eqref{Qformula} with the coefficients $a_0$, $a_{jk}$,  $b_{j,\alpha}$,
 and $c$ being smooth.  Given this, the next step is to observe that Proposition \ref{middlebehavior}
 implies that if we restrict the coefficients to the region where $R\le (\pi-T)/4$, then the $a_0$
 and $a_{jk}$ must vanish to second order at $P_0$ while the $b_{j,\alpha}$ must vanish to first order there.
 If we combine the two steps, we conclude that they have this order of vanishing at $P_0$ when regarded as
 functions of ${\Bbb R}\times S^3$, which completes the proof.

 To conclude this section, we show that, in order to establish \eqref{nullest},
 it suffices to consider the case where the Cauchy data vanishes, that is, $u_0=u_1=0$, $v_0=v_1=0,$ and where $F$ and $G$ vanish for $T$ near 0. To see this, fix $0<\delta\le 1$ and $R_0<\pi$ so that  $\Kob_*\subset \{R<R_0-\delta\}\,.$  On the set $[0,\delta]\times S^3$, the function $\text{dist}(P,P_0)$ is bounded from below, and thus the estimate \eqref{nullest} restricted to this set is implied by the following \begin{multline*}
 \sum_{|\alpha|\le1}\|(\text{dist}^2(P,P_0)\Gamma)^\alpha {\cal Q}\|_{L^2([0,\delta]\times S^3\backslash\Kob_*)}
 \\
 \le C\Bigl(\|u_0\|_{H^2_D}+\|u_1\|_{H^1_D}+ \sum_{|\alpha|\le1}\|(\text{dist}^2(P,P_0)\Gamma)^\alpha F\|_{L^2(\Cylp\backslash\Kob_*)}\Bigr) \\ \times \Bigl(\|v_0\|_{H^2_D}+\|v_1\|_{H^1}+
 \sum_{|\alpha|\le1}\|(\text{dist}^2(P,P_0)\Gamma)^\alpha G\|_{L^2(\Cylp\backslash\Kob_*)}\Bigr).
 \end{multline*} This estimate is established by separately considering the set $R<R_0$, on which it follows by the Minkowski estimate \eqref{1.17}, and the set $R>R_0$, on which it holds by the local estimates of \cite{S1}. 
 Now fix a function $\psi=\psi(T)$ which vanishes near 0, and $\psi=1$ for $T\ge\delta$. It then suffices to establish \eqref{nullest} with $u$ replaced by $\psi u$, which vanishes near $T=0$, and $F$ replaced by $\psi F+\bigl[\square_g,\psi\bigr]u$.  By combining Minkowski energy estimates on $R<R_0$ with $\Cylp$ energy estimates on $R>R_0$, together with the fact that $\psi'$ is supported in $[0,\delta]$,  we obtain the bound \begin{multline*} \sum_{|\alpha|\le1}\|(\text{dist}^2(P,P_0)\Gamma)^\alpha  \bigl(\psi F+\bigl[\square_g,\psi\bigr]u\bigr)\|_{L^2(\Cylp\backslash\Kob_*)} \\ \le C\Bigl(\|u_0\|_{H^2_D}+\|u_1\|_{H^1_D}+ \sum_{|\alpha|\le1}\|(\text{dist}^2(P,P_0)\Gamma)^\alpha F\|_{L^2(\Cylp\backslash\Kob_*)}\Bigr). \end{multline*} Similar comments hold for $v$. This completes the reduction of \eqref{nullest} to the case of vanishing Cauchy data.   \newsection{Unit neighborhoods of the obstacle in Minkowski space}
 
 In this section and the next, we establish the part of \eqref{nullest} (with vanishing Cauchy data) where the norm on the left hand side is taken over a set in the Einstein diamond corresponding to a unit neighborhood of the obstacle in Minkowski space. Precisely, we fix $A>0$ so that $\Kob\subset |x|\le A/2\,.$ Let $$ Y_+={\cal P}\bigl( \{(t,x):t>0\,, x\in\Rt\backslash \Kob\text{ and }|x|<A\} \bigr)\,. $$ We similarly define $Y^*_+$ and $Y^{**}_+$ by replacing $A$ by $2A$ and $3A$ respectively. In a neighborhood of the tip $P_0=(\pi,0)$ of the Einstein diamond, $Y_+$ has the form $$Y_+=\{(T,X)\in \Cylp\backslash \Kob_*: \, R\le A_0(T)\,(\pi-T)^2\}\,,$$ where $A_0(T)$ is smooth and nonvanishing at $P_0\,.$ Also, on $Y_+$, we have $$\Omega\approx (\pi-T)^2\approx (1+t)^{-2}\,.$$ Similar comments hold for $Y^*_+$ and $Y^{**}_+$.  We shall establish estimates on $Y_+$ by pulling them back to Minkowski space. The following will play a key role.   \begin{proposition}\label{prop3.1}
 Let $\tilde \square \tilde u=\tilde F$ and $\tilde \square \tilde v=\tilde G$  be solutions of the Dirichlet-wave equation
 in ${\Bbb R}_+\times {\Bbb R}^3\backslash \Kob$ with  vanishing Cauchy data. Assume that $\tilde F(t,x)=\tilde G(t,x)=0$ when $|x|>2A$.  Then if $N\ge 0$ is fixed,
 \begin{multline}\label{3.1}
 \sum_{|\alpha|\le1}\|(1+t)^ND^\alpha Q(d\tilde u,d\tilde v)\|_{L^2(|x|<2A, t>0)}
 \\
 \le C \sum_{|\alpha|\le1}\|(1+t)^{N/2}D^\alpha \tilde F\|_{L^2(dxdt)} \sum_{|\alpha|\le1}\|(1+t)^{N/2}D^\alpha \tilde G\|_{L^2(dxdt)}.
 \end{multline}
 \end{proposition}
 
 We are using the notation of the preceding section.  In particular
 $\tilde \square = \partial^2/\partial t^2-\Delta$ is the standard D'Alembertian in Minkowski space, and
 $D=(\partial/\partial t,\dots, \partial/\partial x_3)$.
 
 Proposition \ref{prop3.1} follows in a straightforward way from Proposition
 \ref{prop1.3} and Theorem \ref{theorem1.6}.  Before giving the simple proof,
 let us see its relevance for \eqref{nullest}.
 
 We shall apply \eqref{3.1} when $N=2$. We observe that we have the following estimate, if
 $\tilde u=\Penrose^*(\Omega u)$ and $\tilde v=\Penrose^*(\Omega v)$
 \begin{equation}\label{3.2}
 \sum_{|\alpha|\le1}\|(\text{dist}^2(P,P_0)\Gamma)^\alpha {\cal Q}\|_{L^2(Y_+)}
 \le C\sum_{|\alpha|\le1}\|(1+t)^2D^\alpha Q(d\tilde u,d\tilde v)\|_{L^2(|x|<A,t>0)}\,.
 \end{equation}
  This is easy to verify.  One first notes that the usual measure on $\Cyl$ gets pulled back to $\Omega^4dxdt$
 in view of
 \eqref{metrics}. Also, $\text{dist}^2(P,P_0)\approx \Omega$ in $Y_+$.  Therefore,
 Proposition \ref{middlebehavior}, \eqref{conformal}, and \eqref{calQ} imply that the left side of \eqref{3.2}
 is controlled by
 $$\sum_{|\alpha|\le1}\|\Omega^2D^\alpha [\Omega^{-3}Q(d\tilde u,d\tilde v)]\|_{L^2(|x|<A,t>0))}
 \le C\sum_{|\alpha|\le1}\|\Omega^{-1}D^\alpha Q(d\tilde u,d\tilde v)\|_{L^2(|x|<A,t>0))}.$$
 Since $\Omega^{-1}=O((1+t)^2)$ in $Y_+$, this yields our assertion \eqref{3.2}.
 
 To proceed, we first recall that
 $$\tilde \square \tilde u=\tilde F=\Omega^3F.$$
 We cannot apply \eqref{3.1}, though, since $\tilde F$ does not have the required support property.
 To get around this, choose $\eta\in C^\infty_0({\Bbb R}^3)$ satisfying $\eta(x)=1$ if $|x|<A$
 and $\eta(x)=0$ if $|x|>2A$.
 We then let
 $\tilde u^1=\eta \tilde u$,
 and note that
 $$\tilde \square \tilde u^1=2\nabla_x\eta\cdot \nabla_x\tilde u+(\Delta_x\eta)\tilde u+\eta \tilde F=\tilde F^1$$
 has the required support properties. Also, note that $\tilde u^1=\tilde u$ when
 $|x|\le A$.
 
 If we do the same for $\tilde v$, then \eqref{3.1} and \eqref{3.2} imply that
 \begin{multline}\label{3.3}
 \sum_{|\alpha|\le1}\|(\text{dist}^2(P,P_0)\Gamma)^\alpha {\cal Q}(u,du;v,dv)\|_{L^2(Y_+)}
 \\
 \le C\Bigl(\sum_{|\alpha|\le 1}\|(1+t)D^\alpha \tilde F^1\|_{L^2(dxdt)}\Bigr)
 \Bigl(\sum_{|\alpha|\le 1}\|(1+t)D^\alpha \tilde G^1\|_{L^2(dxdt)}\Bigr)
 \\
 \le C \Bigl(\sum_{|\alpha|\le 1}\|(1+t)D^\alpha F\|_{L^2(|x|\le 2A, t>0)}
 +\sum_{|\alpha|\le2}\|(1+t)D^\alpha u\|_{L^2(|x|\le 2A, t>0)}\Bigr)
 \\
 \times \Bigl(\sum_{|\alpha|\le 1}\|(1+t)D^\alpha G\|_{L^2(|x|\le 2A, t>0)}
 +\sum_{|\alpha|\le2}\|(1+t)D^\alpha v\|_{L^2(|x|\le 2A, t>0)}\Bigr).
 \end{multline}
 
 If we repeat the proof of \eqref{3.2}, we conclude that the term
 involving $F$ in the right is controlled by
 %\begin{multline*}
\begin{equation*}
 \sum_{|\alpha|\le1}\|(\pi-T)^{-1}\Omega^{|\alpha|}\Omega^{3}\Omega^{-2}\Gamma^\alpha F\|_{L^2(Y^*_+)}
 %\\
 \approx \sum_{|\alpha|\le1}\|\text{dist}(P,P_0)(\text{dist}^2(P,P_0)\Gamma)^\alpha F\|_{L^2(Y^*_+)}\,.
\end{equation*}
 %\end{multline*}
 Similarly, the term involving $\tilde u$ in the right side of \eqref{3.3} is dominated by
 $$
 \sum_{|\alpha|\le1}\|\text{dist}^{1+2(|\alpha|-2)}(P,P_0)\Gamma^\alpha u\|_{L^2(Y^*_+)}.$$
 Putting all of this together gives us the following
 
 \begin{proposition}\label{prop3.2}
 \begin{multline}\label{3.4}
 \sum_{|\alpha|\le1}\|(\text{dist}^2(P,P_0)\Gamma)^\alpha {\cal Q}\|_{L^2(Y_+)}
 \\
 \le C
 \Bigl(\sum_{|\alpha|\le 1}\|\text{dist}^{1+2|\alpha|}(P,P_0)\Gamma^\alpha F\|_{L^2(Y^*_+)}
 +\sum_{|\alpha|\le2}\|\text{dist}^{1+2(|\alpha|-2)}(P,P_0)\Gamma^\alpha u\|_{L^2(Y^*_+)}\Bigr)
 \\
 \times
 \Bigl(\sum_{|\alpha|\le 1}\|\text{dist}^{1+2|\alpha|}(P,P_0)\Gamma^\alpha G\|_{L^2(Y^*_+)}
 +\sum_{|\alpha|\le2}\|\text{dist}^{1+2(|\alpha|-2)}(P,P_0)\Gamma^\alpha v\|_{L^2(Y^*_+)}\Bigr).
 \end{multline}
 \end{proposition}
 
 \noindent{\bf Proof of Proposition \ref{prop3.1}.} We first note that for
 $j=1,2,3,\dots$ the local estimate \eqref{1.17}
 %from \cite{SS2}
 gives
 \begin{multline}\label{3.5}
 \sum_{|\alpha|\le1}\|(1+t)^ND^\alpha Q(d\tilde u, d\tilde v)\|_{L^2(|x|<2A, t\in [j-1,j])}
 \\
 \le C
 \Bigl((1+j)^{N/2}\bigl[\|\tilde u(j,\cdot)\|_{H^2_D(|x|<2A+1)}+ \|\partial_t \tilde u(j,\cdot)\|_{H^1_D(|x|<2A+1)}\bigr]
 \\
 +\sum_{|\alpha|\le1}\|(1+t)^{N/2}D^\alpha \tilde F\|_{L^2(\{(t,x): t\in [j-1,j]\})}
 \Bigr)
 \\
 \times
 \Bigl((1+j)^{N/2}\bigl[\|\tilde v(j,\cdot)\|_{H^2_D(|x|<2A+1)}+ \|\partial_t \tilde v(j,\cdot)\|_{H^1_D(|x|<2A+1)}\bigr]
 \\
 +\sum_{|\alpha|\le1}\|(1+t)^{N/2}D^\alpha \tilde G\|_{L^2(\{(t,x): t\in [j-1,j]\})}
 \Bigr).
 \end{multline}
 
 On the other hand, the decay estimates \eqref{energydecay},
 yields
 \begin{multline*}
 (1+j)^{N/2}\bigl(\|\tilde u(j,\cdot)\|_{H^2_D(|x|<2A+1)}+ \|\partial_t \tilde u(j,\cdot)\|_{H^1_D(|x|<2A+1)}\bigr)
 \\
 \le C_{N}
 \sum_{|\alpha|\le1}\|e^{-c_0(j-t)}(1+t)^{N/2}D^\alpha \tilde F\|_{L^2(\{(t,x):t<j\})},
 \end{multline*}
 for some $c_0>0$.
 The same estimate works for the other factor in the right side
 of \eqref{3.5}.
 
 If we combine the last two sets of inequalities and square the left side we get
 \begin{multline*}
 \sum_{|\alpha|\le1}\|(1+t)^N D^\alpha Q(d\tilde u, d\tilde v)\|^2_{L^2(|x|<2A, t\in [j-1,j])}
 \\
 \le C_{N}\sum_{1\le i\le j}\sum_{|\alpha|
 \le1}e^{-2c_0(j-i)}\|(1+t)^{N/2} D^\alpha \tilde F\|^2_{L^2([i-1,i]\times {\Bbb R}^3\backslash \Kob)}
 \\
 \times
 \sum_{1\le i\le j}\sum_{|\alpha|
 \le1}e^{-2c_0(j-i)}\|(1+t)^{N/2} D^\alpha \tilde G\|^2_{L^2([i-1,i]\times {\Bbb R}^3\backslash \Kob)}
 \end{multline*}
 If we sum this inequality over $j$, we obtain \eqref{3.1}.
     
 \newsection{Localized energy estimates}
 
 We shall be able to handle the terms in \eqref{3.4} involving $\Gamma^\alpha u$ and $\Gamma^\alpha v$ using the following
 
 \begin{proposition}\label{prop4.1}  If $u$ and $F$ are as in \eqref{linear}, and $u_0=u_1=0,$ then \begin{equation}\label{4.1}
 \sum_{|\alpha|\le2}\|\text{dist}^{1+2(|\alpha|-2)}(P,P_0)\Gamma^\alpha u\|_{L^2(Y^*_+)}
 \le C
 \sum_{|\alpha|\le1}\|(\text{dist}^{2}(P,P_0)\Gamma)^\alpha F\|_{L^2(\Diap)}.
 \end{equation}
 \end{proposition}
 
 If we combine this with Proposition \ref{prop3.2} we conclude that the analog of \eqref{nullest}
 holds if on the left hand side the norm is taken over $Y_+\,.$

 To prove Proposition \ref{prop4.1} we require a couple of Minkowski space estimates
 for solutions of
 \begin{equation}\label{lineartilde}
 \begin{cases}
 \tilde \square \tilde u(t,x)=\tilde F(t,x), \quad (t,x)\in {\Bbb R}_+\times {\Bbb R}^3\backslash \Kob
 \\
 \tilde u(t,x)=0, \quad x\in \partial \Kob
 \\
 \tilde u(0,x)=0, \quad  \partial_t \tilde u(0,x)=0.
 \end{cases}
 \end{equation}
 
 The first estimate we need follows immediately from \eqref{energydecay}.  It says that if $N>0$
 and $A<+\infty$ are fixed then
 \begin{multline}\label{4.2}
 \sum_{|\alpha|\le2}
 \|(1+t)^N D^\alpha \tilde u\|_{L^2(|x|<2A,t>0)}
 \le C \sum_{|\alpha|\le1}\|(1+t)^N D^\alpha \tilde F\|_{L^2(dxdt)},
 \\
 \text{if } \, \tilde F(t,x)=0, \, \, \text{when } \,
 |x|>2A.
 \end{multline}
 
We shall also require the following Minkowski space estimate which is useful when $\tilde F$  vanishes near the boundary.
 
 \begin{proposition}\label{prop4.3}  Suppose that $\tilde u$ solves \eqref{lineartilde}.  Assume further
 that $\tilde F(t,x)=0$ when $|x|<5A/2$.
 Let $\tilde \square \tilde u_f=\tilde F$ be the solution of the corresponding free wave equation. Then given fixed $N>0$
 \begin{equation}\label{4.3}
 \sum_{|\alpha|\le 2}\|(1+t)^ND^\alpha \tilde u\|_{L^2(|x|<2A,t>0)}
 \le C \sum_{|\alpha|\le2}\|(1+t)^ND^\alpha \tilde u_f\|_{L^2(|x|<5A/2,t>0)}.
 \end{equation}
 \end{proposition}
 
 To handle $u_f$ in \eqref{4.3} we shall use the following result.  As we shall see it follows easily from
 Huygen's principle and standard estimates for the free wave equation.
 
 \begin{lemma}\label{lemma4.4}
 Let $\tilde u_f$ solve the free (no obstacle) wave equation $\square u_f=F$ on the Einstein cylinder, with zero Cauchy data.
 Then
 \begin{equation}\label{4.4}
 \sum_{|\alpha|\le2}\|\text{dist}^{1+2(|\alpha|-2)}(P,P_0)\Gamma^\alpha u_f\|_{L^2(Y^{**}_+)} \le C \sum_{|\alpha|\le1}
 \|(\text{dist}^2(P,P_0)\Gamma)^\alpha F\|_{L^2(\Diap)}.
 \end{equation}
 \end{lemma}
 
 Let us momentarily postpone the proofs and see how \eqref{4.2}-\eqref{4.4} can be used to prove
 Proposition \ref{prop4.1}.
 
 For this we first notice that the arguments giving \eqref{3.2} imply that
 $$ \sum_{|\alpha|\le2}\|\text{dist}^{1+2(|\alpha|-2)}(P,P_0)\Gamma^\alpha u
 \|_{L^2(Y^*_+)}
 \le C\sum_{|\alpha|\le2}\|(1+t)D^\alpha \tilde u\|_{L^2(|x|<2A,t>0)}.
 $$ Here $u$ and $\tilde u$ are identified as above.
 
 We then fix $\theta\in C^\infty_0({\Bbb R}^3)$ satisfying $\theta(x)=1$ if $|x|\le 5A/2$ and $\theta(x)=0$ if
 $|x|\ge 3A$.  Using $\theta$ we split
 $$\tilde F=\theta \tilde F+(1-\theta)\tilde F=\tilde G+\tilde H,$$
 We then decompose $\tilde u=\tilde v+\tilde w$, where $\tilde \square \tilde v=\tilde G$. Note then that the forcing term $\tilde H$
 for $\tilde w$ vanishes for $|x|\le 5A/2$.  Also, the preceding inequality yields
 \begin{multline}\label{4.5}
 \sum_{|\alpha|\le2}\|\text{dist}^{1+2(|\alpha|-2)}(P,P_0)\Gamma^\alpha u \|_{L^2(Y^*_+)}
 \\
 \le C\sum_{|\alpha|\le2}\|(1+t)D^\alpha \tilde v\|_{L^2(|x|<2A, t>0)}+
 C\sum_{|\alpha|\le2}\|(1+t) D^\alpha \tilde w\|_{L^2(|x|<2A, t>0)}.
 \end{multline}
 
 If we use \eqref{4.2}
 %Proposition \ref{prop4.2}
 we get that
 \begin{align}\label{4.6}
 \sum_{|\alpha|\le2}\|(1+t)D^\alpha \tilde v\|_{L^2(|x|<2A, t>0)}
 &
 \le C\sum_{|\alpha|\le1} \|(1+t)D^\alpha \tilde G\|_{L^2(dxdt)}
 \notag\\
 &\le C\sum_{|\alpha|\le1} \|(1+t)D^\alpha \tilde F\|_{L^2(|x|<3A,t>0)}
 \notag\\
 &\le C\sum_{|\alpha|\le1} \|\text{dist}^{1+2|\alpha|}(P,P_0)\Gamma^\alpha F\|_{L^2(Y^{**}_+)}. \end{align}
 
 Similarly, if we use Proposition \ref{prop4.3} we get that
 \begin{equation}\label{4.7}
 \sum_{|\alpha|\le2}\|(1+t)D^\alpha \tilde w\|_{L^2(|x|<2A, t>0)}
 \le C\sum_{|\alpha|\le2}\|(1+t)D^\alpha \tilde w_f\|_{L^2(|x|<3A, t>0)},
 \end{equation}
 if $\tilde \square \tilde w_f=\tilde H$ is the solution of the free wave equation.
 If then $w_f$ and $H$ denote the images
 of $\Omega^{-1}\tilde w_f$ and $\Omega^{-3}\tilde H$ in the Einstein diamond, then
 Lemma \ref{lemma4.4} and the above arguments yield
 \begin{align}\label{4.8}
 \sum_{|\alpha|\le2}\|(1+t)D^\alpha \tilde w_f\|_{L^2(|x|<3A,t>0)}&\le C\sum_{|\alpha|\le2}
 \|\text{dist}^{1+2(|\alpha|-2)}(P,P_0)\Gamma^\alpha w_f\|_{L^2(Y^{**}_+)}
 \notag\\
 &\le C \sum_{|\alpha|\le1}\|(\text{dist}^2(P,P_0)\Gamma)^\alpha H\|_{L^2(\Diap)}
 \notag\\
 &\le C \sum_{|\alpha|\le1}\|(\text{dist}^2(P,P_0)\Gamma)^\alpha F\|_{L^2(\Diap)}.
 \end{align}
 In the last step we used the fact that
 $$\sum_{|\alpha|=1}|\Gamma^\alpha H|\le C\sum_{|\alpha|\le1}\text{dist}^{2|\alpha|-2}(P,P_0)
 |\Gamma^\alpha F|.$$
 
Combining \eqref{4.5}-\eqref{4.8} yields Proposition \ref{prop4.1}.  

\noindent{\bf Proof of Proposition \ref{prop4.3}.} Let us first write
 $$\tilde u=\tilde u_f-\tilde u_r,$$
 where $\tilde u_f$ is the solution to the free (i.e., boundaryless) wave equation
 $\tilde \square \tilde u_f=\tilde F$ in ${\Bbb R}^{1+3}_+$, and where $\tilde u_r$ is the reflection term.
 
 To make use of the support assumptions, let us fix $\rho\in C^\infty_0({\Bbb R}^3)$
 satisfying $\rho(x)=1$ if $|x|<2A$ and $\rho(x)=0$ if $|x|>5A/2$.  Then clearly
 \begin{equation}\label{refl}
 \tilde u(t,x)=\rho(x)\tilde u_f(t,x)-\tilde u_r(t,x) \quad \text{if } \, |x|<2A.
 \end{equation}
 
 We next observe that $\rho\tilde u_f -\tilde u_r$ vanishes on ${\Bbb R}_+\times \partial\Kob$, and
 $$\tilde \square (\rho \tilde u_f -\tilde u_r) =\rho\tilde F-2\nabla_x\rho\cdot\nabla_x\tilde u_f
 -(\Delta_x\rho)\tilde u_f=-2\nabla_x\rho\cdot \nabla_x \tilde u_f-(\Delta_x\rho)\tilde u_f,$$ since the support assumptions imply that $\rho\tilde F=0$.
 For the sake of notation, let $\tilde G$ denote the right side of this equation and also set
 $\tilde w=\rho \tilde u_f-\tilde u_r$.  Note that $\tidle G(t,x)=0$ if $|x|>5A/2$.
 
 We now argue as in the proof of Proposition \ref{prop3.1}.  If $j=0,1,\dots$ then \eqref{refl}
 and \eqref{energydecay} yield
 \begin{align*}
 \sum_{|\alpha|\le2}\|D^\alpha \tilde u\|_{L^2(|x|<2A, t\in [j,j+1])}
 &=
 \sum_{|\alpha|\le2}\|D^\alpha \tilde w\|_{L^2(|x|<2A, t\in [j,j+1])}
 \\
 &\le C\,\sum_{|\alpha|\le1}\|e^{-c(j+1-t)}\,D^\alpha \tilde G\|_{L^2(\{(t,x): 0<t<j+1\})}.
 \end{align*}
 This yields \eqref{4.3} since
 $$\sum_{|\alpha|\le1}\|D^\alpha \tilde G(t,\cdot)\|_{L^2(dx)}
 \le C\sum_{|\alpha|\le2}\|D^\alpha \tilde u_f(t,\cdot)\|_{L^2(|x|<5A/2)}.
 $$

 \noindent{\bf Proof of Lemma \ref{lemma4.4}.}
 Write $[0,\pi)=\cup_{j>0} I_j$ where $I_j$ are intervals $[a_j,b_j]$ with $a_{j+1}=b_j$ and
 $$|I_j|\approx (\pi-b_j)^2.$$
 Given $I_j$ let
 \begin{equation}\label{4.15}
 \Lambda_j=\{(T,X)\in \Diap: \, T+R\in I_j\},
 \end{equation}
 where, as in Section 2, $R$ is the distance from $X$ to the north pole.  Then if $T\in I_j$, by Huygen's principle, the energy inequality, and the fact that the $\Gamma_j$ commute with the D'Alembertian, there is a uniform constant $B$, depending on $A$, so that
 $$
 \sum_{|\alpha|=2}\|\Gamma^\alpha u_f\|_{L^\infty L^2(Y^{**}_+:T\in I_j)}
 \le C\sum_{|k-j|<B}\sum_{|\alpha|=1}\|\Gamma^\alpha F\|_{L^1L^2(\Lambda_k)}.
 $$
 By the Schwarz inequality we can dominate the first term on the right by
 \begin{multline*} 
\sum_{|k-j|<B}\sum_{|\alpha|=1}
\|\text{dist}^2(P,P_0)\Gamma^\alpha F\|_{L^2(\Lambda_k)}\times
 \Bigl(\int_0^T(\pi-s)^{-4}\, ds\Bigr)^{1/2}
 \\
 \le
 C\sum_{|k-j|<B}\sum_{|\alpha|=1}(\pi-T)^{-3/2}\|\text{dist}^2(P,P_0)\Gamma^\alpha F\|_{L^2(\Lambda_k)}.
 \end{multline*} 
Using these inequalities, we conclude that
\begin{multline*}
\sum_{|\alpha|=2}\|\text{dist}(P,P_0)
%&
\Gamma^\alpha u_f\|_{L^2(Y^{**}_+:T\in I_j)}
%\\
%&
\le C(\pi-T)^2\sum_{|\alpha|=2}\|\Gamma^\alpha  u_f\|_{L^\infty L^2(Y^{**}_+:T\in I_j)}
\\
%&
\le  C\sum_{|k-j|<B}\sum_{|\alpha|=1}\|\text{dist}^{5/2}(P,P_0)\Gamma^\alpha
 F\|_{L^2(\Lambda_k)}
\,.
\end{multline*}
 From this we deduce that
 $$
 \sum_{|\alpha|=2}\|\text{dist}(P,P_0)\Gamma^\alpha u_f\|_{L^2(Y^{**}_+)}
 \le C\sum_{|\alpha|=1}\|\text{dist}^{5/2}(P,P_0)\Gamma^\alpha F\|_{L^2(\Diap)} $$ which shows that the terms on the left hand side of \eqref{4.4} corresponding to $|\alpha|=2$ satisfy stronger estimates than those
 asserted by the Lemma.
 
 The same bounds hold for the terms involving $\Gamma^\alpha u_f$, $|\alpha|=1$.  To see this, we first use H\"older's
 inequality to see that if $T$ is fixed, then
 $$\|\text{dist}^{-1}(P,P_0)\Gamma^\alpha u_f(T,\cdot)\|_{L^2(X:(X,T)\in Y^{**}_+)}
 \le C\|\text{dist}(P,P_0)\Gamma^\alpha u_f(T,\cdot)\|_{L^6(X:(X,T)\in Y^{**}_+)}.$$
 Next, by the energy inequality and Sobolev embedding, if $(\square_g+1)w=H$ with vanishing Cauchy data then
 $$\|w(T,\cdot)\|_{L^6(S^3)}\le  C\int_0^T\|H(S,\cdot)\|_{L^2(S^3)}\, dS.$$
 Consequently, if we use Huygen's principle and repeat our earlier arguments we find that
 $$
 \sum_{|\alpha|=1}\|\Gamma^\alpha u_f\|_{L^\infty L^6(Y^{**}_+:T\in I_j)}
 \le C\sum_{|k-j|<B}\sum_{|\alpha|=1}\|\Gamma^\alpha F\|_{L^1L^2(\Lambda_k)}
 \,. $$ Therefore, we can dominate the terms in \eqref{4.4} with $|\alpha|=1$ by the same bounds as for $|\alpha|=2$.
 
 It remains to handle the terms with $|\alpha|=0\,.$ We may bound $$\|\text{dist}(P,P_0)^{-3}u_f\|^2_{L^2(Y^{**}_+:T\in I_j)} \le C\,\|u_f\|^2_{L^\infty_TL^6_X(Y^{**}_+:T\in I_j)} \le C\sum_{|k-j|<B}\|F\|^2_{L^2(\Lambda_k)}\,. $$ Summing over $j$ yields the desired bound, completing the proof of Lemma \ref{lemma4.4}.

 \newsection{Complement of image of unit neighborhood of obstacle}
 
 To handle the complement of $Y_+$ in the Einstein cylinder, we shall use the following estimates for the free wave equation. 
 \begin{proposition}\label{prop5.0} Let $(\square_g+1)u_f=F$, $(\square_g+1)v_f=G$ be solutions of the
 free (no obstacle) wave equation in the Einstein cylinder $\Cylp$. As before, let
 $${\cal Q}(u_f,du_f;v_f,dv_f)=\Omega^{-3}Q(d\Omega u_f, d\Omega v_f).$$ Then
 \begin{multline*}
 \sum_{|\alpha|\le1}\|(\text{dist}^2(P,P_0)\Gamma)^\alpha {\cal Q}\|_{L^2(\Cylp)}
 \\
 \le C\Bigl(\|u_f(0,\cdot)\|_{H^2(S^3)}+\|\partial_Tu_f(0,\cdot)\|_{H^1(S^3)}+ \sum_{|\alpha|\le1}\|(\text{dist}^2(P,P_0)\Gamma)^\alpha F\|_{L^1L^2(\Cylp)}\Bigr) \\ \times \Bigl(\|v_f(0,\cdot)\|_{H^2(S^3)}+\|\partial_Tv_f(0,\cdot)\|_{H^1(S^3)}+
 \sum_{|\alpha|\le1}\|(\text{dist}^2(P,P_0)\Gamma)^\alpha G\|_{L^1L^2(\Cylp)}\Bigr).
 \end{multline*}
 \end{proposition}
  The first step in establishing Proposition \ref{prop5.0} is to observe that it suffices to consider the case where the Cauchy data of $u_f$ and $v_f$ vanish; that is, $u_f(0,\cdot)=\partial_Tu_f(0,\cdot)=0$, and similarly for $v_f$. This follows by a similar (but simpler) reduction to that at the end of Section 2.  Consequently, we are reduced to establishing the following estimate, in the case of vanishing Cauchy data, 
\begin{multline}\label{5.1}
 \sum_{|\alpha|\le1}\|(\text{dist}^2(P,P_0)\Gamma)^\alpha {\cal Q}\|_{L^2(\Cylp)}
 \\
 \le C \sum_{|\alpha|\le1}\|(\text{dist}^2(P,P_0)\Gamma)^\alpha F\|_{L^1L^2(\Cylp)} \sum_{|\alpha|\le1}\|(\text{dist}^2(P,P_0)\Gamma)^\alpha G\|_{L^1L^2(\Cylp)}.
 \end{multline}  We postpone the proof of the estimate \eqref{5.1} and state a few consequences. For simplicity, we state the next three estimates in the case of vanishing Cauchy data. The first just follows from the
 Schwarz inequality,
 \begin{multline}\label{5.2}
 \sum_{|\alpha|\le1}\|(\text{dist}^2(P,P_0)\Gamma)^\alpha {\cal Q}\|_{L^2(\Cylp)}
 \\
 \le C\sum_{|\alpha|\le1}\|(\text{dist}^2(P,P_0)\Gamma)^\alpha F\|_{L^2(\Cylp)}
 \sum_{|\alpha|\le1}\|(\text{dist}^2(P,P_0)\Gamma)^\alpha G\|_{L^2(\Cylp)}.
 \end{multline}
 
 As in the preceding section, we can obtain better weighted estimates when
 the forcing terms are supported near the boundary.  One such estimate that we shall need is the following
 \begin{multline}\label{5.3}
 \sum_{|\alpha|\le1}\|(\text{dist}^2(P,P_0)\Gamma)^\alpha {\cal Q}\|_{L^2(\Cylp)}
 \\
 \le C\sum_{|\alpha|\le1}\|\text{dist}^{1+2|\alpha|}(P,P_0)\Gamma^\alpha F\|_{L^2(\Cylp)}
 \sum_{|\alpha|\le1}\|\text{dist}^{2|\alpha|}(P,P_0)\Gamma^\alpha G\|_{L^2(\Cylp)},
 \\
 \text{if }\,\text{support}(F)\subseteq Y_+\,.
 \end{multline}
 Similarly,
 \begin{multline}\label{5.4}
 \sum_{|\alpha|\le1}\|(\text{dist}^2(P,P_0)\Gamma)^\alpha {\cal Q}\|_{L^2(\Cylp)}
 \\
 \le C\sum_{|\alpha|\le1}\|\text{dist}^{1+2|\alpha|}(P,P_0)\Gamma^\alpha F\|_{L^2(\Cylp)}
 \sum_{|\alpha|\le1}\|\text{dist}^{1+2|\alpha|}(P,P_0)\Gamma^\alpha G\|_{L^2(\Cylp)},
 \\
 \text{if }\,\text{support}(F,G)\subseteq Y_+\,.
 \end{multline}
 
 The proof of \eqref{5.3} and \eqref{5.4} uses a decomposition that is similar to the one
 employed in the proof of Lemma \ref{lemma4.4}. Let $I_j$ be as in the proof of Lemma \ref{lemma4.4}, and let $\Psi_j$ be a partition of unity such that $\Psi_j=1$ on $I_j$ and $\Psi_j$ is supported in the doubled interval $I_j^*$. Also, let
 $\Lambda^+_j$ be the subset of $\Cylp$ where $T-R\in I_j$.  Then if $F$ is supported in $Y_+$,
 we can use \eqref{5.1} and Huygen's principle to see that
 \begin{multline*}
 \sum_{|\alpha|\le1}\|(\text{dist}^2(P,P_0)\Gamma)^\alpha {\cal Q}\|_{L^2(\Lambda^+_j)}
 \\
 \le C\sum_{|k-j|<B}
 \sum_{|\alpha|\le1}\|(\text{dist}^2(P,P_0)\Gamma)^\alpha (\Psi_kF)\|_{L^1L^2(\Cylp)}
 \sum_{|\alpha|\le1}\|(\text{dist}^2(P,P_0)\Gamma)^\alpha G\|_{L^1L^2(\Cylp)},
 \end{multline*}
 if $B$ is a large fixed constant as before.
 Since $\Gamma^\alpha\Psi_k=O\bigl(\text{dist}^{-2|\alpha|}(P,P_0)\bigr)$ on the support of $F$, we see from
 the Schwarz inequality that
 $$
 \sum_{|\alpha|\le1}\|(\text{dist}^2(P,P_0)\Gamma)^\alpha (\Psi_kF)\|_{L^1L^2(\Cylp)}
 \le C\sum_{|\alpha|\le1}\|\text{dist}^{1+2|\alpha|}(P,P_0)\Gamma^\alpha F\|_{L^2(\Cylp:T\in I^*_k)},
 $$
 assuming, as above, that $F$ is supported in $Y_+$. By combining the last two inequalities and applying
 the Schwarz inequality, we get
 \begin{multline*}
 \sum_{|\alpha|\le1}\|(\text{dist}^2(P,P_0)\Gamma)^\alpha {\cal Q}\|_{L^2(\Lambda^+_j)}
 \\
 \le C\sum_{|k-j|<B}\sum_{|\alpha|\le1}\|\text{dist}^{1+2|\alpha|}(P,P_0)\Gamma^\alpha F\|_{L^2(\Cylp:T\in I_k^*)}
 \sum_{|\alpha|\le1}\|\text{dist}^{2|\alpha|}\Gamma^\alpha G\|_{L^2(\Cylp)},
 \\
 \text{if }\, \text{support}(F)\subseteq Y_+.
 \end{multline*}
 This implies \eqref{5.3}.
 The proof of \eqref{5.4} is similar.
 
 We now show that Proposition \ref{prop5.0}
 implies the inequality \eqref{nullest} where the norm on the left hand side is taken over  $\Cylp\backslash Y_+.$
 Together with the results of Sections 3 and 4, this will complete the proof of the inequality \eqref{nullest}.  For this, let $\beta$ be the pushforward to $\Cylp$ of the function $1-\eta(2x)$. Thus,
 $$
 \begin{cases}
 \beta=1 \quad \text{on } \Cylp\backslash Y_+
 \\
 \beta=0\quad \text{on a neighborhood of }\Kob^*
 \\
 \Gamma^\alpha \beta = O\bigl(\text{dist}^{-2|\alpha|}(P,P_0)\bigr).
 \end{cases}
 $$
 We set
 $$u_f=\beta u, \quad v_f=\beta v,$$
 so that $u_f=u$ and $v_f=v$ on $Y_+^c.$ Additionally,
 they solve the free wave equations
 \begin{align*} (\square_g+1)u_f & =\beta F+[\square_g,\beta]u = F_0+F_1\\ (\square_g+1)v_f & =\beta G+[\square_g,\beta]v = G_0+G_1 \end{align*} with vanishing Cauchy data.
 We then can write $u_f=u_{f,0}+u_{f,1}$ and $v_f=v_{f,0}+v_{f,1}$, where the pieces solve the free wave equations
 $$(\square_g+1)u_{f,j}=F_j, \quad (\square_g+1)v_{f,j}=G_j, \quad j=0,1.$$
 
 Note also that
 \begin{equation}\label{5.5}
 \sum_{|\alpha|\le1}\|\text{dist}^{2|\alpha|}(P,P_0)\Gamma^\alpha F_0\|_{L^2(\Cylp)}
 \le C\sum_{|\alpha|\le1}\|\text{dist}^{2|\alpha|}(P,P_0)\Gamma^\alpha F\|_{L^2(\Cylp)}.
 \end{equation}
 Furthermore, since $F_1$ vanishes on $Y_+^c$, and
 $$\sum_{|\alpha|\le1}|(\text{dist}^2(P,P_0)\Gamma)^\alpha [\square_g,\beta]u|
 \le C\sum_{|\alpha|\le2}\text{dist}^{2(|\alpha|-2)}(P,P_0)|\Gamma^\alpha u|
 $$
 Proposition \ref{prop4.1} yields
 \begin{multline}\label{5.6}
 \sum_{|\alpha|\le1}\|\text{dist}^{1+2|\alpha|}(P,P_0)\Gamma^\alpha F_1\|_{L^2(\Cylp)}
 \le C\sum_{|\alpha|\le2}
 \|\text{dist}^{1+2(|\alpha|-2)}(P,P_0)\Gamma^\alpha u\|_{L^2(Y_+)}
 \\
 \le C\sum_{|\alpha|\le1}\|\text{dist}^{2|\alpha|}(P,P_0)\Gamma^\alpha F\|_{L^2(\Cylp)}.
 \end{multline}
 To proceed, note that \eqref{5.2} and \eqref{5.5} yield
 \begin{multline*}
 \sum_{|\alpha|\le1}\|(\text{dist}^2(P,P_0)\Gamma)^\alpha {\cal
 Q}(u_{f,0},du_{f,0};v_{f,0},dv_{f,0})\|_{L^2(\Cylp)}
 \\
 \le
 C\sum_{|\alpha|\le1}\|\text{dist}^{2|\alpha|}(P,P_0)\Gamma^\alpha F\|_{L^2(\Cylp)}
 \sum_{|\alpha|\le1}\|\text{dist}^{2|\alpha|}(P,P_0)\Gamma^\alpha G\|_{L^2(\Cylp)}.
 \end{multline*}
 Similarly, \eqref{5.3}, \eqref{5.5} and \eqref{5.6} give
 \begin{align*}
 \sum_{|\alpha|\le1}&\|(\text{dist}^2(P,P_0)\Gamma)^\alpha {\cal Q}(u_{f,1},du_{f,1};v_{f,0},dv_{f,0})\|_{L^2(\Cylp)}
 \\
 &+\sum_{|\alpha\le1} \|(\text{dist}^2(P,P_0)\Gamma)^\alpha {\cal Q}(u_{f,0},du_{f,0};v_{f,1},dv_{f,1})\|_{L^2(\Cylp)}
 \\
 &\le C\sum_{|\alpha|\le1}\|\text{dist}^{2|\alpha|}(P,P_0)\Gamma^\alpha F\|_{L^2(\Cylp)} \sum_{|\alpha|\le1}\|\text{dist}^{2|\alpha|}(P,P_0)\Gamma^\alpha G\|_{L^2(\Cylp)}.
 \end{align*}
 Finally, \eqref{5.4} and \eqref{5.6} imply that
 \begin{multline*}
 \sum_{|\alpha|\le1}\|(\text{dist}^2(P,P_0)\Gamma)^\alpha {\cal Q}(u_{f,1},du_{f,1};v_{f,1},dv_{f,1})\|_{L^2(\Cylp)}
 \\
 \le
 C\sum_{|\alpha|\le1}\|\text{dist}^{2|\alpha|}(P,P_0)\Gamma^\alpha F\|_{L^2(\Cylp)}
 \sum_{|\alpha|\le1}\|\text{dist}^{2|\alpha|}(P,P_0)\Gamma^\alpha G\|_{L^2(\Cylp)}.
 \end{multline*}
 
 If we combine the last three inequalities we conclude that the analog of \eqref{nullest} holds if the norm on
 the left is taken over $Y_+^c$.

 %\bigskip
 \noindent {\bf Proof of the estimate \eqref{5.1}}
  In what follows, we simplify our notation somewhat by letting
 $$|h'|=\sum_{|\alpha|=1}|\Gamma^\alpha h|, \quad \text{and } \, |h''|=\sum_{|\alpha|=2}|\Gamma^\alpha h|.$$
 In view of Proposition \ref{nulltransformationprop},  Proposition \ref{prop5.0} follows as a consequence of the following three lemmas.

 \begin{lemma}\label{lemma5.2}  Suppose that $Q$ is a standard null form on the cylinder and that
 $u_f$ and $v_f$ are solutions of standard inhomogeneous wave equations $(\square_g +1)u_f=F$ and
 $(\square_g+1)v_f=G$ with zero initial data. Then \begin{multline}\label{5.7}
 \sum_{|\alpha|\le1}\|\text{dist}^2(P,P_0)(\text{dist}^2(P,P_0)\Gamma)^\alpha Q(du_f,dv_f)\|_{L^2(\Cylp)}
 \\
 \le
 C\sum_{|\alpha|\le1}
 \|\text{dist}^{2|\alpha|}(P,P_0)\Gamma^\alpha F\|_{L^1L^2(\Cylp)}
 \|\text{dist}^{2|\alpha|}(P,P_0)\Gamma^\alpha G\|_{L^1L^2(\Cylp)}.
 \end{multline}
 \end{lemma}
 
 \begin{lemma}\label{lemma5.3}  Let $u_f$ and $v_f$ be as above then
 \begin{multline}\label{5.8}
 \|\text{dist}^3(P,P_0)u_fv_f''\|_{L^2(\Cylp)}
 +\|\text{dist}^3(P,P_0)u_f'v_f'\|_{L^2(\Cylp)}
 \\
 \le C
 \sum_{|\alpha|\le1}\bigl\|\text{dist}^{2|\alpha|}(P,P_0)\Gamma^\alpha F\bigr\|_{L^1L^2(\Cylp)}
 \|\sum_{|\alpha|\le1}\bigl\|\text{dist}^{2|\alpha|}(P,P_0)\Gamma^\alpha G\bigr\|_{L^1L^2(\Cylp)}.
 \end{multline}
 \end{lemma}
 
 \begin{lemma}\label{lemma5.4}  If $u_f$ and $v_f$ are as above
 \begin{multline}\label{5.9}
 \|\text{dist}(P,P_0)u_f'v_f\|_{L^2(\Cylp)} + \|u_fv_f\|_{L^2(\Cylp)}
 \\
 \le C
 \sum_{|\alpha|\le1}\bigl\|\text{dist}^{2|\alpha|}(P,P_0)\Gamma^\alpha F\bigr\|_{L^1L^2(\Cylp)}
 \sum_{|\alpha|\le1}\bigl\|\text{dist}^{2|\alpha|}(P,P_0)\Gamma^\alpha G\bigr\|_{L^1L^2(\Cylp)}.
 \end{multline}
 \end{lemma}
 
 \noindent {\bf Proof of Lemma \ref{lemma5.2}}
 We apply an estimate from \cite{S1}, which says that
 $$\sum_{|\alpha|\le1}\|\Gamma^\alpha Q(du_f,dv_f)\|_{L^2(\Cylp)}
 \le C\sum_{|\alpha|\le1}\|\Gamma^\alpha F\|_{L^1L^2(\Cylp)}\sum_{|\alpha|\le1}\|\Gamma^\alpha G\|_{L^1L^2(\Cylp)}.
 $$  Next, we fix a partition of unity $\sum \beta(2^js)=1$, $s>0$
 with $\text{supp }\beta\subset [1/2,2]$, and
 let $F_j=\beta(2^j\text{dist}(P,P_0))F$.
 From the preceeding estimate, Huygen's principle, and the Schwarz inequality, one sees that
 for some fixed $B$, the following holds for $k=0,1,2,\dots$
 \begin{align*}
 \sum_{|\alpha|=1}\|&\text{dist}^4(P,P_0)\Gamma^\alpha
 Q\|_{L^2(\text{dist}(P,P_0)\approx 2^{-k})}
 \\
 &\le C \sum_{|\alpha|\le1}\sum_{j\le k+B}\|2^{-2k}\Gamma^\alpha F_j\|_{L^1L^2}
 \sum_{|\alpha|\le1}\sum_{j\le k+B}\|2^{-2k}\Gamma^\alpha G_j\|_{L^1L^2}
 \\
 &\le C \sum_{|\alpha|\le1}\sum_{j\le
 k+B}2^{-2(k-j)}\|\text{dist}^{2|\alpha|}(P,P_0)\Gamma^\alpha  F\|_{L^1L^2(\text{dist}(P,P_0)\approx 2^{-j})}
 \\
 &\qquad\qquad\qquad \times
 \sum_{|\alpha|\le1}\sum_{j\le k+B}2^{-2(k-j)}\|\text{dist}^{2|\alpha|}(P,P_0) \Gamma^\alpha G\|_{L^1L^2(\text{dist}(P,P_0)\approx 2^{-j})}.
 \end{align*}
 This yields the estimate for the $|\alpha|=1$ terms on the left side of \eqref{5.7}.
 
 To estimate the term with $\alpha=0$ on the left side of  \eqref{5.7}, we 
 let $u_{f,k}$ be the solution of the inhomogeneous wave equation $(\square_g+1)u_{f,k}=F_k$.  We then can write
 $$Q=\sum_{j\le l}Q(du_{f,j},dv_{f,l})+\sum_{j>l}Q(du_{f,k},dv_{f,l}).$$
 
 Since the two terms are similar, we shall only estimate the first one.  We use the  following estimate from \cite{S1}.
 \begin{equation}
 \label{5.10}
 \|Q\|_{L^2(\Cylp)}\le C\sum_{|\alpha|\le1}\|\Gamma^\alpha F\|_{L^1L^2(\Cylp)}\|G\|_{L^1L^2(\Cylp)},
 \end{equation}
 We apply this to obtain
 \begin{multline*}
 \bigl\|\text{dist}^2(P,P_0)\sum_{j\le l}Q(du_{f,j},dv_{f,l})\bigr\|_{L^2(\text{dist}(P,P_0)\approx
 2^{-k})}
 %\\
 \le C\sum_{|\alpha|\le1}\sum_{j\le l\le k+B}2^{-2k}\|\Gamma^\alpha F_j\|_{L^1L^2}
 \|G_l\|_{L^1L^2}
 \\
 \le C\sum_{|\alpha|\le1}\sum_{j\le l\le k+B}2^{-2(k-j)}
 \|\text{dist}^{2|\alpha|}(P,P_0)\Gamma^\alpha F\|_{L^1L^2(\text{dist}(P,P_0)\approx 2^{-j})}
 \|G\|_{L^1L^2(\text{dist}(P,P_0)\approx 2^{-l})} \\ \le C\sum_{|\alpha|\le1}\sum_{j\le k+B}2^{-(k-j)}
 \|\text{dist}^{2|\alpha|}(P,P_0)\Gamma^\alpha F\|_{L^1L^2(\text{dist}(P,P_0)\approx 2^{-j})}
\\
\times
\sum_{l\le k+B}2^{-(k-l)}\|G\|_{L^1L^2(\text{dist}(P,P_0)\approx 2^{-l})}
 \end{multline*}
 This yields the desired estimate for the $|\alpha|=0$ term in the left of \eqref{5.7},  completing the proof of Lemma \ref{lemma5.2}.
 
 \medskip
 
 In the proof of Lemma \ref{lemma5.3}, we shall need to make use of the following estimate.
 
 \begin{lemma}\label{lemma5.5}  Let $B_\varepsilon$ denote a spherical cap of radius $\varepsilon>0$ in $S^3$.  Then
 $$\|h\|_{L^\infty(B_{\varepsilon})}\le C\varepsilon^{1/2}\|h'\|_{L^6(B_{\varepsilon})}
 + C\varepsilon^{-1/2}\|h\|_{L^6(B_{\varepsilon})}.$$
 \end{lemma}
 
 To prove this one first notices that it follows from the Euclidean version, which in turn
 follows from the $\varepsilon=1$ case and a simple scaling argument.
 
 \noindent {\bf Proof of Lemma \ref{lemma5.3}.} We start by estimating the
 first term on the left hand side of \eqref{5.8} since it is the more difficult. We use the preceeding lemma to obtain \begin{align*} \|u_f\|_{L^\infty(\text{dist}(P,P_0)\approx 2^{-k})} & \le  \sum_{j\le k+B}\|u_{f,j}\|_{L^\infty(\text{dist}(P,P_0)\le 2^{-k})} \\ &\le\sum_{j\le k+B} 2^{-k/2}\|u_{f,j}'\|_{L^\infty L^6}+2^{k/2}\|u_{f,j}\|_{L^\infty L^6}\\ &\le\sum_{j\le k+B} 2^{-k/2}\|F_j'\|_{L^1L^2}+2^{k/2}\|F_j\|_{L^1L^2}\\ &\le\sum_{j\le k+B} 2^{-k/2}\|F\|_{L^1L^2(\text{dist}(P,P_0)\approx 2^{-j})}+2^{k/2}\|F\|_{L^1L^2(\text{dist}(P,P_0)\approx 2^{-j})} \end{align*} As a result, for fixed
$k=0,1,2,\dots$,
\begin{align*}
\|&\text{dist}^3(P,P_0)u_fv_f''\|_{L^2(\text{dist}(P,P_0)\approx 2^{-k})}
\\
&\le C\,2^{-3k-k/2}\|v_f''\|_{L^\infty L^2(\text{dist}(P,P_0)\approx 2^{-k})} \|u_f\|_{L^\infty(\text{dist}(P,P_0)\approx 2^{-k})}
\\
&\le C\,2^{-2k} \sum_{j\le k+B}\|G'\|_{L^1L^2(\text{dist}(P,P_0)\approx
2^{-j})}
\\
&\quad\times\Bigl(2^{-2k}\sum_{j\le k+B}\|F'\|_{L^1L^2(\text{dist}(P,P_0)\approx
2^{-j})}+2^{-k}\sum_{j\le k+B}\|F\|_{L^1L^2(\text{dist}(P,P_0)\approx 2^{-j})}\Bigr)
\\
&\le C \sum_{j\le k+B}2^{-2(k-j)}\|\text{dist}^{2}(P,P_0)G'\|_{L^1L^2(\text{dist}(P,P_0)\approx
2^{-j})}
\\
&\quad
\times \sum_{j\le k+B}
\Bigl(
2^{-2(k-j)}\|\text{dist}^{2}(P,P_0)F'\|_{L^1L^2(\text{dist}(P,P_0)\approx2^{-j})}
%\\
%&\qquad\qquad\qquad\qquad
+
%\sum_{j\le k+B}
2^{-k}\|F\|_{L^1L^2(\text{dist}(P,P_0)\approx 2^{-j})}\Bigr) .
\end{align*}
 This implies that the first term on the left side of \eqref{5.8} satisfies the stated estimate.
 
 For the second term in \eqref{5.8}, we first use H\"older's inequality to obtain
 \begin{multline*}
 \|\text{dist}^3(P,P_0)u_f'v_f'\|_{L^2(\text{dist}(P,P_0)\approx 2^{-k})}
 \\
 \le C2^{-4k}\|u_f'\|_{L^\infty L^6(\text{dist}(P,P_0)\approx 2^{-k})}
 \|v_f'\|_{L^\infty L^6(\text{dist}(P,P_0)\approx 2^{-k})}
 \end{multline*}
 As above, we can bound 
\begin{multline*} 2^{-2k}\|u_f'\|_{L^\infty L^6(\text{dist}(P,P_0)\approx 2^{-k})}\\ \le 
 C\sum_{j\le k+B}
\bigl(
2^{-2(k-j)}\|\text{dist}^{2}(P,P_0)F'\|_{L^1L^2(\text{dist}(P,P_0)\approx 2^{-j})} +2^{-k}\|F\|_{L^1L^2(\text{dist}(P,P_0)\approx 2^{-j})}
\bigr)
 \end{multline*}  
which leads to the desired bounds for the
 remaining term, completing the proof of Lemma \ref{lemma5.3}.
 
 \medskip
 
 \noindent{\bf Proof of Lemma \ref{lemma5.4}}  Let $u_{f,j}$ and $v_{f,l}$ be as in the proof of
 Lemma \ref{lemma5.2}.  We then can write
 $$u_f'v_f=\sum_{j\le l}u'_{f,j}v_{f,l}+\sum_{j>l}u'_{f,j}v_{f,l}.$$
 To handle the terms with $j\le l$, we note that H\"older's inequality and the above arguments yield
 \begin{align*}
 \sum_{j\le l\le k+B}
 &\|\text{dist}(P,P_0)u'_{f,j}
 v_{f,l}\|_{L^2(\text{dist}(P,P_0)\approx 2^{-k})}
 \\
 &\le C\,2^{-2k}\sum_{j\le l\le k+B} \|u'_{f,j}\|_{L^\infty L^6} \|v_{f,l}\|_{L^\infty L^6}
 \\
 &\le C \sum_{j\le l\le k+B}
 \Bigl(2^{-2(k-j)}\|\text{dist}^{2}(P,P_0)F'\|_{L^1L^2(\text{dist}(P,P_0)\approx
 2^{-j})} +2^{-k}\|F\|_{L^1L^2}\Bigr)
 \\
 &\qquad\qquad\qquad\qquad\times
  \|G\|_{L^1L^2(\text{dist}(P,P_0)\approx 2^{-l})}.
 \end{align*}
 Since $k-j\ge k-l$, this implies the estimate for these terms. 
 To handle the terms where $j>l$ we write
 \begin{align*}
 \sum_{l\le j\le k+B}&\|\text{dist}(P,P_0)u'_{f,j}v_{f,l}\|_{L^2(\text{dist}(P,P_0)\approx 2^{-k})}
\\
 &\le C\,2^{-3k/2}\sum_{l\le j\le k+B}\|u'_{f,j}\|_{L^\infty L^2(\text{dist}(P,P_0)\approx 2^{-k})} \|v_{f,l}\|_{L^\infty(\text{dist}(P,P_0)\le 2^{k})}
 \\
 &\le C\sum_{l\le j\le k+B}\|F\|_{L^1L^2(\text{dist}(P,P_0)\approx 2^{-j})}
\\ &\qquad\qquad\qquad\times
 \Bigl(2^{-2(k-l)}\|\text{dist}(P,P_0)^2G'\|_{L^1L^2(\text{dist}(P,P_0)\approx 2^{-l})}+
 2^{-k}\|G\|_{L^1L^2}\Bigr).
 \end{align*}
 This establishes the desired estimate for the terms with $j>l$, which along with the preceding estimate shows that the first term on the left side of \eqref{5.9}
 satisfies the desired bounds.
 
 To handle the second term on the left hand side of \eqref{5.9}, we apply H\"older's inequality to deduce that
 \begin{align*}
 \|u_fv_f\|_{L^2(\text{dist}(P,P_0)\approx 2^{-k})}
 &\le C\,2^{-k}\|u_f\|_{L^\infty L^6(\text{dist}(P,P_0)\approx 2^{-k})}
 \|v_f\|_{L^\infty L^6(\text{dist}(P,P_0)\approx 2^{-k})}
 \\
 &\le C\, 2^{-k}\|F\|_{L^1L^2(\text{dist}(P,P_0)\ge 2^{-k+B})}\|G\|_{L^1L^2(\text{dist}(P,P_0)\ge 2^{-k+B})}
 \end{align*}
 which implies the desired estimate, completing the proof of Lemma \ref{lemma5.4}.
 
 \newsection{End of Proof of Theorem \ref{theorem2.3}: Weighted Pecher Estimates}
 
 We still need to prove \eqref{weightedenergy}, which is the weighted Pecher inequality
 \begin{multline}\label{6.1}
 \sum_{|\alpha|\le 1}
 \|(\text{dist}^2(P,P_0)\Gamma)^\alpha u\|_{L^8(\Cylp\backslash \Kob_*)}\\
 \le C\Bigl(\|u_0\|_{H^2_D}+\|u_1\|_{H^1_D}+\sum_{|\alpha|\le1}\|(\text{dist}^2(P,P_0)\Gamma)^\alpha
 F\|_{L^2(\Cylp\backslash \Kob_*)}\Bigr).
 \end{multline}
  The local version of the Pecher inequalities for variable coefficient wave equations was established in \cite{MSS}. In particular, that result implies the following result for the free wave equation on the Einstein cylinder \begin{equation}\label{localpecher} \|u_f\|_{L^8(\Cylp)} \le C \bigl(\|u_f(0,\cdot)\|_{H^1}+\|\partial_T u_f(0,\cdot)\|_{L^2}+ \|F\|_{L^1L^2(\Cylp)}\bigr).
 \end{equation} Also, the local version of the Pecher inequalities for the obstacle problem in Minkowski space were established by the authors in \cite{SS1}, $$ \|\tilde u\|_{L^8([0,1]\times\Rt\backslash\Kob)} \le C\bigl(\|\tilde u_0\|_{H^1_D}+\|\tilde u_1\|_{L^2}+ \|\tilde F\|_{L^1L^2([0,1]\times\Rt\backslash\Kob)}\bigr). $$ We also need the version with first order derivatives of $\tilde u$, which follows from the above and an integration by parts argument in $t$, 
\begin{equation}\label{minkpecher} \sum_{|\alpha|\le 1}\|D^\alpha\tilde u\|_{L^8([0,1]\times\Rt\backslash\Kob)} 
%\\ 
\le C\bigl(\|\tilde u_0\|_{H^2_D}+\|\tilde u_1\|_{H^1_D}+ \sum_{|\alpha|\le 1}\|D^\alpha \tilde F\|_{L^1L^2([0,1]\times\Rt\backslash\Kob)}\bigr). 
\end{equation} 
Together, \eqref{localpecher} and \eqref{minkpecher} allow one to reduce the proof of the estimate \eqref{6.1} to the case that $u_0=u_1=0$, following the arguments at the end of the second section of this paper.  The proof of \eqref{6.1} now follows very closely the proof of \eqref{nullest}. As in the proof of that estimate, the first step is to control the norm over $Y_+.$ To do this, we note that $$ \sum_{|\alpha|\le 1} \|(\text{dist}^2(P,P_0)\Gamma)^\alpha u\|_{L^8(Y_+)} \approx \sum_{|\alpha|\le 1} \|(1+t)D^\alpha \tilde u\|_{L^8(|x|\le A,t>0)}. $$ Next, if $\square\tilde u^1=\tilde F^1$, with vanishing Cauchy data, and $\tilde F^1$ is supported in the set $|x|\le 2A$, then energy decay and the local estimates \eqref{minkpecher} imply the following $$ \sum_{|\alpha|\le 1} \|(1+t)D^\alpha \tilde u^1\|_{L^8(|x|\le A,t>0)} \le \sum_{|\alpha|\le 1} \|(1+t)D^\alpha \tilde F^1\|_{L^2(dxdt)}. $$ If we set $\tilde u^1=\eta \tilde u$, then (see \eqref{3.3}, \eqref{3.4}, and Proposition \ref{prop4.1}) we obtain the estimate \begin{multline*}
 \sum_{|\alpha|\le 1}
 \|(\text{dist}^2(P,P_0)\Gamma)^\alpha u\|_{L^8(Y_+)}\\
 \le C\Bigl(\|u_0\|_{H^2_D}+\|u_1\|_{H^1_D}+\sum_{|\alpha|\le1}\|(\text{dist}^2(P,P_0)\Gamma)^\alpha
 F\|_{L^2(\Cylp\backslash \Kob_*)}\Bigr).
 \end{multline*}  To handle the norm over the complement of $Y_+$, we will use the following estimates for the free wave equation on $\Cylp$, \begin{equation}\label{6.4}
 \sum_{|\alpha|\le1}\|(\text{dist}^2(P,P_0)\Gamma)^\alpha u_f\|_{L^8(\Cylp)}
 \le C\sum_{|\alpha|\le1}\|\text{dist}^{2|\alpha|}(P,P_0)\Gamma^\alpha F\|_{L^2(\Cylp)} \end{equation} and the improved estimate for data supported near the boundary \begin{multline}\label{6.5}
 \sum_{|\alpha|\le1}\|(\text{dist}^2(P,P_0)\Gamma)^\alpha u_f\|_{L^8(\Cylp)}
 \le C\sum_{|\alpha|\le 1}\|\text{dist}^{1+2|\alpha|}(P,P_0)\Gamma^\alpha F\|_{L^2(\Cylp)}
 \\ \text{if }\,\text{support}(F)\subseteq Y_+\,.
 \end{multline} The estimate \eqref{6.5} is a consequence of \eqref{6.4} by the same steps as \eqref{5.3} follows from \eqref{5.1}. And, by letting $u_f=\beta u$, the following is a consequence of \eqref{6.4} and \eqref{6.5} \begin{multline*}
 \sum_{|\alpha|\le 1}
 \|(\text{dist}^2(P,P_0)\Gamma)^\alpha u\|_{L^8(\Cylp\backslash Y_+)}\\
 \le C\Bigl(\|u_0\|_{H^2_D}+\|u_1\|_{H^1_D}+\sum_{|\alpha|\le1}\|(\text{dist}^2(P,P_0)\Gamma)^\alpha
 F\|_{L^2(\Cylp\backslash \Kob_*)}\Bigr).
 \end{multline*}  It thus remains only to establish the estimate \eqref{6.4}. For the $|\alpha|=0$ terms, this is just the estimate \eqref{localpecher}. Next, from the fact that $\Gamma^\alpha$ commutes with $\square_g$, we obtain the following $$
 \sum_{|\alpha|=1}\|\Gamma^\alpha u_f\|_{L^8(\text{dist}(P,P_0)\approx 2^{-j})}
 \le C\sum_{|\alpha|=1}\|\Gamma^\alpha F\|_{L^2(\text{dist}(P,P_0)\ge 2^{-j+B})} $$ for $B$ fixed as before. Summing over $j$ yields \eqref{6.4}.  \qed
 %\bigskip
 
 We conclude this section with a simple corollary of our weighted Pecher estimate \eqref{6.1}.
 We first see that, if $u$ is as above, then
 \begin{equation*}
 (\pi-T)|u(T,\cdot)|\le
 C\Bigl(\|u_0\|_{H^2_D}+\|u_1\|_{H^1_D}+\sum_{|\alpha|\le1}\|(\text{dist}^2(P,P_0)\Gamma)^\alpha
 F\|_{L^2(\Cylp\backslash \Kob_*)}\Bigr).
 \end{equation*}
 To prove this, one uses the fact that, if $B_\varepsilon$ is a ball of radius
 $\varepsilon>0$, then
 $$
 \|h\|_{L^\infty(B_\varepsilon)}\le
 C\varepsilon^{1/2}\sum_{|\alpha|=1}\|\Gamma^\alpha h\|_{L^8(B_{\varepsilon})}
 +C\varepsilon^{-1/2}\|h\|_{L^8(B_{\varepsilon})},$$
 which follows from Euclidean estimates.  Given $(T,X)\in\Diap\backslash\Kob_*$, we take $\varepsilon=(\pi-T)^2$, and $B_\epsilon$ to be a ball of radius $\epsilon$ such that $(T,B_\epsilon)$ is contained in $\Diap\backslash\Kob_*$. We then obtain from \eqref{6.1} the following inequality,
 \begin{align*}
 |u(T,X)|
 &\le C\,(\pi-T)\sum_{|\alpha|=1}\|\Gamma^\alpha u\|_{L^8(B_\epsilon)}
 +C\,(\pi-T)^{-1}\|u\|_{L^8(B_\epsilon)}
 \\
 &\le C(\pi-T)^{-1}
 \sum_{|\alpha|\le1}\|\text{dist}^{2|\alpha|}(P,P_0) \Gamma^\alpha u\|_{L^8(\Diap\backslash\Kob_*)}
 \\
 &\le C(\pi-T)^{-1}\Bigl(\|u_0\|_{H^2_D}+\|u_1\|_{H^1_D}+ \sum_{|\alpha|\le1}\|(\text{dist}^2(P,P_0)\Gamma)^\alpha
 F\|_{L^2(\Cylp\backslash \Kob_*)}\Bigr),
 \end{align*}
 as claimed.
 
 We now conclude that the solution to equation \eqref{1.3} decays like
 $1/t$.  For this we note that, if $u$ and $\tilde u={\cal P}^*(\Omega u)$ are identified as before, and $\Penrose(t,x)=(T,X)$, then
 \begin{multline}\label{6.6}
 |\tilde u(t,x)|=|\Omega u(T,X)|\\
 \le C t^{-1}\Bigl(\|u_0\|_{H^2_D}+\|u_1\|_{H^1_D}
 +\sum_{|\alpha|\le1}\|(\text{dist}^2(P,P_0)\Gamma)^\alpha
 F\|_{L^2(\Cylp\backslash \Kob_*)}\Bigr).
 \end{multline}
 This inequality uses the fact that in $\Diap\backslash \Kob_*$,
 we have $|\Omega/(\pi-T)|\le C/t$.
 
 \newsection{Global Existence of $H^2$ Solutions: Proof of Theorem
 \ref{theorem1.1}}
 
 Recall that we need to show that if the Cauchy data satisfies
 \begin{equation}\label{7.1}
 \|f\|_{H^{2,1}_D({\Bbb R}^3\backslash \Kob)}+\|g\|_{H^{1,2}_D({\Bbb
 R}^3\backslash \Kob)}< \varepsilon_0,
 \end{equation}
 with $\varepsilon_0>0$ small, then the equation
 \begin{equation}\label{7.2}
 \begin{cases}
 \square u(t,x)=Q(du(t,x),du(t,x)), \quad (t,x)\in {\Bbb R}_+\times {\Bbb
 R}^3\backslash \Kob
 \\
 u|_{\partial \Kob}=0
 \\
 u|_{t=0}=f, \quad \partial_tu|_{t=0}=g
 \end{cases}
 \end{equation}
 has a unique global solution verifying \eqref{1.9} and \eqref{1.10}.  Note that
 implicit in \eqref{7.1} is that the data satisfy the $H^2$ compatibility
 conditions that both $f$ and $g$ vanish on $\partial \Kob$.
 
 To avoid cumbersome notation we are switching our notation from the last
 several sections.  In this section and the one to follow we do not denote
 functions and derivatives on Minkowski space with a tilde.
 
 The uniqueness assertion follows immediately from \eqref{1.17}.  As we shall
 see, the existence assertion follows easily from Theorem \ref{theorem2.3}.  Precisely, we shall use Theorem \ref{theorem2.3} to solve the corresponding equation on the Einstein cylinder minus the obstacle. Restricting this solution to the Einstein diamond yields a solution to \eqref{7.2}, after pulling back via the Penrose transform. 
 Thus, let
 \begin{equation*}
 u=\Omega {\cal P}^*v
 \end{equation*}
 denote $\Omega$ times the pullback of $v$ via the Penrose map.  (In our
 previous
 notation $u$ would be $\tilde v$.)  Then, as noted before, \eqref{7.2} is
 implied by the following
 \begin{equation}\label{7.3}
 \begin{cases}
 (\square_g+1)v(T,X)={\cal Q}(v(T,X),dv(T,X); v(T,X),dv(T,X)), \quad (T,X)\in
 \Cylp\backslash \Kob_*
 \\
 v(T,X)=0, \quad (T,X)\in \Kob_*
 \\
 v|_{T=0}=f_e, \quad \partial_Tv|_{T=0}=g_e,
 \end{cases}
 \end{equation}
 assuming that $f=\Omega {\cal P}_0^*f_e$ and 
$g=\Omega^2{\cal P}_0^* g_e$, with  ${\cal P}_0$ denoting the restriction of the Penrose map to $t=0$.
 
 To construct a solution of \eqref{7.3} on the Einstein cylinder, we let $v$ denote the solution
 of the following linear equation
 \begin{equation*}
 \begin{cases}
 (\square_g+1)v=F,  \quad (T,X)\in \Cylp\backslash \Kob_*
 \\
 v|_{\partial \Kob_*}=0
 \\
 v|_{T=0}=f_e, \quad \partial_Tv|_{T=0}=g_e,
 \end{cases}
 \end{equation*}
 The existence of a solution to this linear equation on the Einstein diamond  minus the obstacle is obtained from the corresponding solution on Minkowski space. That solution is then 
easily extended to the Einstein cylinder minus the obstacle.  Now let 
$$
 {\cal T} F={\cal Q}(v,dv; v,dv)\,,
 $$ 
where $v$ solves the above linear equation. Finding a solution to \eqref{7.3} is 
thus reduced to finding a fixed point
 for the operator $\cal T$ on the set of $F$ such that
 $$
 \sum_{|\alpha|\le1}\|(\text{dist}^2(P,P_0)\Gamma)^\alpha
 F\|_{L^2(\Cylp\backslash \Kob_*)}\le c_0\,.
 $$
 We next observe that, by Theorem \ref{theorem2.3}, 
 $$
 \sum_{|\alpha|\le1}\|(\text{dist}^2(P,P_0)\Gamma)^\alpha
 {\cal T}(0)
 \|_{L^2(\Cylp\backslash \Kob_*)}\le C\,\epsilon_0\,,
 $$
 and
 \begin{multline*}
 \sum_{|\alpha|\le1}\|(\text{dist}^2(P,P_0)\Gamma)^\alpha
 ({\cal T}F_1-{\cal T}F_2)
 \|_{L^2(\Cylp\backslash \Kob_*)}\\
 \le
 %\sum_{|\alpha|\le1}\|(\text{dist}^2(P,P_0)\Gamma)^\alpha
 %{\cal Q}(v_1+v_2,d(v_1+v_2);v_1-v_2,d(v_1-v_2))
 %\|_{L^2(\Cylp\backslash \Kob_*)}
 \sum_{|\alpha|\le1}\|(\text{dist}^2(P,P_0)\Gamma)^\alpha
 {\cal Q}(v_1-v_2,d(v_1-v_2);v_1,dv_2)
 \|_{L^2(\Cylp\backslash \Kob_*)}
\\
\quad 
 +\sum_{|\alpha|\le1}\|(\text{dist}^2(P,P_0)\Gamma)^\alpha
 {\cal Q}(v_2,dv_2;v_1-v_2,d(v_1-v_2))
 \|_{L^2(\Cylp\backslash \Kob_*)}
\\
 \le
 2\,C\,(c_0+\epsilon_0)
 \cdot
 \sum_{|\alpha|\le1}\|(\text{dist}^2(P,P_0)\Gamma)^\alpha
 (F_1-F_2)\|_{L^2(\Cylp\backslash \Kob_*)}
\,.
 \end{multline*}
 Taking $c_0=2\,C\,\epsilon_0$, then, for $\epsilon_0>0$ small enough, the contraction
 principle yields a fixed point for $\cal T\,,$ hence a solution
 to \eqref{7.3}.
 The corresponding function $u$ defined as above by $u=\Omega \Penrose_*v$
 is then a solution of \eqref{7.2}, and it must verify \eqref{1.9}
 by appealing to \eqref{1.17}, and also
 must satisfy the decay estimate \eqref{1.10} because of \eqref{6.6}.

 \newsection{Global Existence of Smooth Solutions: Proof of Theorem 
 \ref{theorem1.3}}
 To establish Theorem \ref{theorem1.3}, we will show that, if the data $f$ and $g$ are smooth and satisfy the appropriate compatibility conditions to infinite order, then the solution
 $u$ given by Theorem \ref{theorem1.1} belongs to 
 $C^\infty({\Bbb R}_+\times {\Bbb R}^3\backslash \Kob)\,.$ The proof
 is based on the following Lemma.
 
 \begin{lemma}\label{lemma8.1}
 Let $u\in C^j([0,T]; H^{2-j}_D(\Rt\backslash\Kob))\,,j=0,1,2\,,$ be
 a solution to \eqref{1.3}. Then if
 $v\in C^j([0,T]; H^{1-j}_D(\Rt\backslash\Kob))\,,\,j=0,1\,,$ solves
 \begin{equation*}
 \begin{cases}
 \square v=Q(du,dv)+Q(dv,du)+F, \quad (t,x)\in {\Bbb R}_+
 \times \Rt\backslash \Kob
 \\
 v(t,\cdot)|_{\partial\Kob}=0
 \\
 v(0,\cdot)\in H^2_D(\Rt\backslash\Kob), \, \, 
 \partial_tv(0,\cdot)\in H^1_D(\Rt\backslash\Kob),
 \end{cases}
 \end{equation*}
 where $F,DF\in L^1([0,T];L^2(\Rt\backslash\Kob))\,,$ and the
 following a priori assumption holds,
 $$ Q(du,dv),\, Q(dv,du)\in L^1([0,T];L^2(\Rt\backslash\Kob))\,,$$ then
 $v\in C^j([0,T]; H^{2-j}_D(\Rt\backslash\Kob))\,,j=0,1,2\,,$
 and
 $$ DQ(du,dv),\, DQ(dv,du)\in L^1([0,T];L^2(\Rt\backslash\Kob))\,. $$
 \end{lemma}
 {\bf Proof.}
 We use the following estimates from \cite{SS2}, for solutions $u,v$
 to the Cauchy-Dirichlet problem,
 \begin{multline}\label{8.2}
 \|DQ(du,dv)\|_{L^2([0,T]\times\Rt\backslash\Kob)}  \\ \le
 C(u,T)\,
 \Bigl(\|v_0\|_{H^2_D(\Rt\backslash\Kob)}+
 \|v_1\|_{H^1_D(\Rt\backslash\Kob)}
 +\sum_{|\alpha|\le 1}\int_0^T\|D^\alpha\square v(s,\cdot\,)
 \|_{L^2(\Rt\backslash\Kob)}\,ds\,\Bigr),
 \end{multline}
 
 \begin{multline}\label{8.3}
 \|Q(du,dv)\|_{L^2([0,T]\times\Rt\backslash\Kob)}  \\ \le 
 C(u,T)\,\Bigl(
 \|v_0\|_{H^1_D(\Rt\backslash\Kob)}+\|v_1\|_{L^2(\Rt\backslash\Kob)}
 +\int_0^T\|\square v(s,\cdot\,)
 \|_{L^2(\Rt\backslash\Kob)}\,ds\,\Bigr)
\,,
 \end{multline}
 where
 $$
 C(u,T)=
 C\times\Bigl(
 \|u_0\|_{H^2_D(\Rt\backslash\Kob)}+\|u_1\|_{H^1_D(\Rt\backslash\Kob)}
 +\sum_{|\alpha|\le 1}\int_0^T\|D^\alpha\square u(s,\cdot\,)
 \|_{L^2(\Rt\backslash\Kob)}\,ds\,\Bigr)\,.
 $$
 By \eqref{8.2}, if $T'<(2C(u,T))^{-2}\,,$ and $\tilde v$ denotes the
 free solution to the Dirichlet-Cauchy problem $\square \tilde v=0$ with data $(\tilde v(0,\cdot),\partial_t \tilde v(0,\cdot))=(v(0,\cdot),\partial_t v(0,\cdot))\in  
 H^2_D\times H^1_D\,,$
 then the Schwarz inequality shows that the map
 $$
 G\rightarrow Q(du,d(\square^{-1}G+\tilde v))+Q(d(\square^{-1}G+\tilde v),du)+F
 $$
 is a contraction in the norm 
 $\sum_{|\alpha|\le 1}\|D^\alpha G\|_{L^2([0,T']\times\Rt\backslash\Kob)}\,.$ Here, for shorthand, $\square^{-1}G$ denotes the solution of the inhomogeneous wave equation with  forcing term $G$ and zero initial data.
 If $G$ is the unique fixed point of this map, then $v=\square^{-1}G+\tilde v$
 satisfies $\square v=Q(du,dv)+Q(dv,du)+F\,.$ A similar proof using \eqref{8.3}
 shows that there is a unique solution for $0\le t\le T'$
 among $v$ such that
 $\square v\in L^1([0,T'];L^2(\Rt\backslash\Kob))\,.$ By the assumptions
 of the lemma, these solutions coincide, which establishes the result
 for $T$ replaced by $T'$. By energy estimates, 
 $v(T',\cdot\,)\in H^2_D(\Rt\backslash\Kob)\,,$ and
 $\partial_t v(T',\cdot\,)\in H^2_D(\Rt\backslash\Kob)\,.$  The above argument may be repeated to establish the regularity result for $0\le t\le 2T'$, 
 and iteration yields the result for $0\le t\le T$. 
 
{\bf Proof of Theorem \ref{theorem1.3}.}
 Let $u$ denote the solution given by Theorem \ref{theorem1.1}. We shall show that, under the additional assumptions of Theorem \ref{theorem1.3},  if $0<T<\infty$ is fixed then $u\in C^\infty([0,T]\times {\Bbb R}^3\backslash \Kob)$. 
 By truncating the data $f,g$ and using finite propagation velocity,
 we may assume that the data is compactly supported (the constant
 $\epsilon_0$ is of course independent of the support.)
 Let $\psi_j(x)=\partial_t^j u(0,x)$ 
 denote the function obtained by formally differentiating $u$ at
 $t=0$ and expressing the result in terms of $f$ and $g$. By assumption,
 $\psi_j\in C_0^\infty(\Rt\backslash\Kob)$, and $\psi_j$ vanishes
 on $\partial\Kob\,.$ 
 Then, for finite $T$, 
 $u\in C^j([0,T]; H^{2-j}_D(\Rt\backslash\Kob))\,,j=0,1,2\,.$
 We differentiate equation \eqref{1.3} in $t$ to obtain the following
 equation for $v=\partial_t u\,,$
 $$
 \begin{cases}
 \square v=Q(du,dv)+Q(dv,du), \quad (t,x)\in {\Bbb R}_+\times \Rt\backslash \Kob
 \\
 v(t,\cdot)|_{\partial\Kob}=0
 \\
 v(0,\cdot)\in H^2_D(\Rt\backslash\Kob), \, \, 
 \partial_tv(0,\cdot)\in H^1_D(\Rt\backslash\Kob).
 \end{cases}
 $$
 By Lemma \ref{lemma8.1},
 $v\in C^j([0,T]; H^{2-j}_D(\Rt\backslash\Kob))\,,j=0,1,2\,,$
 and $DQ(du,dv),\,DQ(dv,du)\in L^2([0,T]\times\Rt\backslash\Kob)\,.$
 
 Suppose that we have shown
 $\partial_t^k u\in C^j([0,T]; H^{2-j}_D(\Rt\backslash\Kob))\,,j=0,1,2\,,$
 for $0\le k<m\,.$
 We differentiate \eqref{1.3} $m$ times in $t$ to obtain the following
 equation for $v=\partial_t^m u\,,$
 $$
 \begin{cases}
 \square v=Q(du,dv)+Q(dv,du)+F, \quad (t,x)\in {\Bbb R}_+\times \Rt\backslash \Kob
 \\
 v(t,\cdot)|_{\partial\Kob}=0
 \\
 v(0,\cdot)\in H^2_D(\Rt\backslash\Kob), \, \, 
 \partial_tv(0,\cdot)\in H^1_D(\Rt\backslash\Kob)\,,
 \end{cases}
 $$
 where 
 $$
 F=\sum_{0<j<m}{m\choose j}\bigl(Q(d\partial_t^j u,d\partial_t^{m-j}u) +Q(d\partial_t^{m-j}u,d\partial_t^j u) \bigr).
 $$
 By the induction step, and \eqref{8.2} and \eqref{8.3}, we have that
 $F,DF\in L^2([0,T]\times\Rt\backslash\Kob)\,,$ and
 $Q(du,dv)\,,Q(dv,du)\in L^2([0,T]\times\Rt\backslash\Kob)\,.$ Consequently,
 $\partial_t^m u\in C^j([0,T]; H^{2-j}_D(\Rt\backslash\Kob))\,,$ for $j=0,1,2\,.$
 It follows that
 $$
 u\in C^\infty([0,T];H^2_D(\Rt\backslash\Kob))\,,\quad
 Q(du,du)\in C^\infty([0,T];H^1(\Rt\backslash\Kob))\,,
 $$
 where, as before, $H^k({\Bbb R}^3\backslash \Kob)$ denotes Sobolev space of restrictions of elements of $H^k(\Rt)\,.$
 
 We next obtain spatial regularity for $u$ from the equation \eqref{1.3}.
 We have
 $$
 \begin{cases}
 \Delta u =\partial_t^2 u+Q(du,du)\in
 C^\infty([0,T];H^1(\Rt\backslash\Kob))\\
 u(t,\cdot)|_{\partial\Kob}=0\,.
 \end{cases}
 $$
 By elliptic regularity, thus $u\in C^\infty([0,T];H^3_D(\Rt\backslash\Kob))\,.$
 Since $H^k(\Rt\backslash\Kob)$ is an algebra under pointwise
 multiplication for $k\ge 2$, it follows that 
 $Q(du,du)\in C^\infty([0,T];H^2(\Rt\backslash\Kob))\,.$
 A simple induction now establishes that 
 $u\in C^\infty([0,T];H^k_D(\Rt\backslash\Kob))$ for all $k\,.$
 
 {\bf Remark.} If it is assumed that the initial data 
 $\psi_j=\partial_t^j u(0,\cdot\,)$ belongs to $H^2_D(\Rt\backslash\Kob)$
 for all $j$, (no smallness assumption is necessary on the norm,
 except for $j=0,1$), then it is not necessary to truncate the data to apply
 Lemma \ref{lemma8.1}, and we may conclude that the solution of
 Theorem \ref{theorem1.1} satisfies
 $u\in C^\infty([0,T];H^k_D(\Rt\backslash\Kob))$ for all $k$.

 \end{document}